\documentclass[11pt]{article}
\usepackage{amsfonts,latexsym,amsmath,amscd,geometry}
\usepackage{amssymb}
\usepackage{latexsym}
\usepackage{color}
\usepackage{url}
\usepackage{hyperref}

\newcommand \nc{\newcommand}
\newtheorem{theorem}{Theorem}[section]
\newtheorem{lemma}[theorem]{Lemma}

\newtheorem{remark}[theorem]{Remark}
\newtheorem{problem}[theorem]{Problem}

\nc{\ba}{\begin{array}}\nc{\ea}{\end{array}}
\nc{\be}{\begin{eqnarray}}\nc{\ee}{\end{eqnarray}}
\nc{\beq}{\begin{equation}}\nc{\eeq}{\end{equation}}
\nc{\bex}{\begin{eqnarray*}}\nc{\eex}{\end{eqnarray*}}
\nc{\btm}{\begin{theorem}} \nc{\etm}{\end{theorem}}
\nc{\blm}{\begin{lemma}} \nc{\elm}{\end{lemma}}

\nc{\R}{\mathbb{R}} 
\nc{\ve}{\varepsilon}
\def\grad{\nabla}
\nc\Tr{\mathrm{tr}}
\DeclareMathOperator{\diver}{div}

\def\cJ{\mathcal J}
\def\RR{\mathbb R}
\def\bfA{{\mathbf A}}
\def\bfB{{\mathbf B}}
\def\bfQ{{\mathbf Q}}
\def\bfR{{\mathbf R}}
\def\bfS{{\mathbf S}}
\def\bfG{{\mathbf G}}

\def\bfv{{\mathbf v}}
\def\obfA{\overline{\mathbf A}}
\def\obfQ{\overline{\mathbf Q}}
\makeatletter
\DeclareFontFamily{OMX}{MnSymbolE}{}
\DeclareSymbolFont{MnLargeSymbols}{OMX}{MnSymbolE}{m}{n}
\SetSymbolFont{MnLargeSymbols}{bold}{OMX}{MnSymbolE}{b}{n}
\DeclareFontShape{OMX}{MnSymbolE}{m}{n}{
    <-6>  MnSymbolE5
   <6-7>  MnSymbolE6
   <7-8>  MnSymbolE7
   <8-9>  MnSymbolE8
   <9-10> MnSymbolE9
  <10-12> MnSymbolE10
  <12->   MnSymbolE12
}{}
\DeclareFontShape{OMX}{MnSymbolE}{b}{n}{
    <-6>  MnSymbolE-Bold5
   <6-7>  MnSymbolE-Bold6
   <7-8>  MnSymbolE-Bold7
   <8-9>  MnSymbolE-Bold8
   <9-10> MnSymbolE-Bold9
  <10-12> MnSymbolE-Bold10
  <12->   MnSymbolE-Bold12
}{}

\let\llangle\@undefined
\let\rrangle\@undefined
\DeclareMathDelimiter{\llangle}{\mathopen}%
                     {MnLargeSymbols}{'164}{MnLargeSymbols}{'164}
\DeclareMathDelimiter{\rrangle}{\mathclose}%
                     {MnLargeSymbols}{'171}{MnLargeSymbols}{'171}
\makeatother
\begin{document}

\title{Mean Field Control and Mean Field Game Models \\ with Several Populations}

\author{Alain Bensoussan\\
\small International Center for Decision and Risk Analysis\\
\small Jindal School of Management, University of Texas at Dallas\thanks{axb046100@utdallas.edu ; also with the College of Science and Engineering, Systems Engineering
and Engineering Management, City University Hong Kong. Research supported
by the National Science Foundation under grant DMS-1612880 and grant
from the SAR Hong Kong RGC GRF 11303316. } \\
Tao Huang\\
\small Department of Mathematics, Wayne State University, Detroit, MI, USA\thanks{gq3481@wayne.edu}\\
Mathieu Lauri{\`e}re\\
\small Department of Operations Research and Financial Engineering, \\ 
\small Princeton
University, Princeton, NJ, USA\thanks{lauriere@princeton.edu ; partially supported by NSF \#DMS- 1716673 and ARO \#W911NF-17-1-0578.}}
%\maketitle

\date{}

\maketitle

\begin{abstract}

In this paper, we investigate the interaction of two populations with a large number of indistinguishable agents. 
The problem consists in two levels: the interaction between agents of a same population, and the interaction between the two populations.
In the spirit of mean field type control (MFC) problems and mean field games (MFG), each population is approximated by a continuum of infinitesimal agents. We define four different problems in a general context and interpret them in the framework of MFC or MFG. By calculus of variations, we derive formally in each case the adjoint equations for the necessary conditions of optimality. Importantly, we find that in the case of a competition between two coalitions, one needs to rely on a system of master equations in order to describe the equilibrium. Examples are provided, in particular linear-quadratic models for which we obtain systems of ODEs that can be related to Riccati equations.
\end{abstract}

%\tableofcontents

\section{Introduction}
\setcounter{equation}{0}
\setcounter{theorem}{0}

% subsec-subsec-subsec-subsec-subsec-subsec-subsec-subsec-subsec %
% subsec-subsec-subsec-subsec-subsec-subsec-subsec-subsec-subsec %
\subsection{General introduction}

The evolution of a large group of interacting agents, who are trying to realize a certain goal either from an individual or a social viewpoint, is an important and interesting question in mathematics, physics and many other fields. When the number of agents grows to infinity, it becomes extremely hard to keep track of all the agent-to-agent interactions and to study the resulting global behavior. However, by assuming that every agent has the same importance, the impact of each single agent's choice on the group decreases as the size of the population increases. So, in order to efficiently approximate the global evolution of the group, one can replace the influences of all the other players on a given agent by their average influence. This is called \textit{mean field approach}, whose name is borrowed from statistical physics. This approach is valid under some assumptions (in particular the one asserting that all the players have indistinguishable roles among the population) and allows to replace the microscopic viewpoint by a macroscopic one. The main advantage of such an approximation is that the macroscopic description is more tractable and in particular amenable to numerical treatment. Moreover, the larger the number of individuals, the more accurate the approximation.

%There are two important theories in mean field approach, 
Based on the idea of mean field approach, two important theories have recently emerged: \textit{mean field games} (MFG) and \textit{mean field (type) control} problems (MFC). Both of them study the behavior of a typical agent \textit{via} the evolution of her state and the cost induced by this evolution.
One can characterize the optimal control of this representative player by two coupled equations: a forward one, describing her dynamics (or, equivalently, the dynamics of the population), and a backward one, describing the evolution of her value function.
Depending on whether one uses analytical or stochastic methods, partial differential equations (PDEs) or stochastic differential equations (SDEs) are obtained respectively. For more details, the reader can refer to~\cite{MR3134900} for the general theory, \cite{MR3559742} for the regularity theory of the PDE system, and~\cite{MR3752669,MR3753660} for the stochastic analysis.

There are several important differences between these two theories that have to be emphasized.

\smallskip
(1) Mean field games correspond to the limit of differential games, in which one wants to find Nash equilibria, when the number of players tend to infinity.
They describe the global behavior of the group resulting from the selfish choices that individuals are making so as to minimize a certain cost (or maximize a certain benefit). This cost depends on the state of the agent and also on the statistical distribution of all the agents' states, that is, the global state of the system. 
Mean field game models have been introduced on the one hand by Lasry and Lions~\cite{MR2269875,MR2271747,MR2295621,PLL-CDF,MR2762362}  (see also~\cite{Cardaliaguet-2013-notes} for a written account of Lions' lectures at Coll{\`e}ge de France), and on the other hand by Caines, Huang and Malham{\'e}~\cite{HuangCainesMalhame-2003-individual-mass-wireless,HuangCainesMalhame-2004-large-LQG,MR2346927,HuangMalhameCaines-2006-NCE-physics,MR2344101,MR2352434}. They have since then attracted a lot of interest. Some important applications can be found in
	finance~\cite{MR2295621,MR2784835,MR3500451,MR3325272},
	economy~\cite{Gueant-2013-phdthesis,MR3363751,MR2647032,MR3359708,MR3268061,MR3268060},
	systemic risk~\cite{MR3032938,MR3325083},
 	energy production~\cite{MR2762362,MR3683681}, 
	or crowd motion~\cite{PLL-CDF,LachapelleWolfram-2011-MFG-congestion-aversion,YAJML}.

(2) Another kind of asymptotic regime leads to mean field type control problems, which are stochastic optimal control problems where the cost function and the parameters of the dynamics depend on the law of the controlled stochastic process. In general, such problems correspond to the control of a large number of agents by a global planner. The problem is to find an optimal feedback rule that she can provide to all the agents, who then implement it in a distributed fashion. This type of problem can also be regarded as the problem of a single player who tries to optimize a cost involving the law of her own state, which evolves with a Mc-Kean Vlasov (MKV) dynamics. Depending on which of the two viewpoints one chooses, this theory is referred to \textit{mean field (type) control} (e.g., by Bensoussan, Frehse, and Yam) or \textit{control of Mc-Kean Vlasov dynamics} (e.g., by Carmona and Delarue).
It has found applications such as risk management, portfolio management, or cybersecurity~\cite{MR3405978,MR3451175,Pfeiffer-2016-risk-merton,MR3575619}.

\medskip
In this paper, we will investigate the interaction of two populations with a large number of indistinguishable agents, which is a natural extension of the aforementioned theories. In this type of problems, two (or more) large populations interact and the behavior of each group is approximated using a mean field approach. This idea has been broached by several authors since the original contribution of Huang, Caines and Malham{\'e}~\cite{MR2346927}. 
In \cite{MR3127148}, Feleqi has derived an adjoint system of Hamilton-Jacobi-Bellman and Kolmogorov-Fokker-Plank for the ergodic mean field game theory of several populations, by letting the number of members of each group go to infinity.
In~\cite{MR3333058,MR3660463} and~\cite{CirantBardi-Uniqueness-SeveralPop-arxiv}, existence and uniqueness results for this type of systems with Neumann boundary conditions have been proved, in the stationary and the dynamic case respectively.
In \cite{MR3597009}, Achdou, Bardi and Cirant introduced a MFG model describing the interactions between two populations in urban settlements and residential choice. They showed the existence of solutions for both stationary and evolution cases with periodic boundary conditions, and provided some numerical simulations. 
For a synthetic presentation, we refer the reader to the monographs~\cite{MR3134900} (Chapter 8) and~\cite{MR3752669} (Chapter 7).  To the best of our knowledge, mean field control problems with several populations have been considered only in~\cite{LachapelleWolfram-2011-MFG-congestion-aversion} and~\cite{MR3763083}, where the authors introduced models for crowd motion (with local and non-local interactions respectively). They studied optimality conditions and provided numerical results.

Despite the increasing research activities on this topic, a global viewpoint was still missing: multi-population MFC problems have been considered only on some examples, and even in the more studied multi-population MFG setting, the problems were of a relatively special type due to the form of the Hamiltonians.
 The goal of the present paper is to introduce a general framework to tackle multi-population mean field control problems and mean field games. 
To this end, we consider two types of interactions, cooperation or competition, between two populations and between agents of one population. At a heuristic level, we can summarize as follows the different cases.

 (1) In the same population, if the agents cooperate each other, one obtains a mean field control problem (or control of McKean-Vlasov dynamics), whereas if they compete, one obtains a mean field game. In this paper, we only consider the case that different groups have the same type of interaction between their own agents.
 
 (2) For two populations, if the interaction is cooperation, there is only one single objective function to optimize and the equilibrium between the two groups should be of social type. If it is competition, each group optimizes its own objective function and we look for a Nash equilibrium between the two populations.

Therefore, we will study four different cases of two population interactions in this paper. For each of them, we will describe the control problems and derive the associated system of PDEs by calculus of variation. Moreover the adjoint equations can also be deduced from the so-called \textit{master equation}~\cite{MR3258261,MR3501391,MR3343705,MR3652408,MR3332710,MR3575613,MR3631380}. We will also explain generalizations of this approach for two populations.

For the sake of clarity, we will consider the interactions between only two populations but the ideas could be generalized to a larger number of populations. 
To alleviate the notations, we will assume that each population represents the same proportion of the total population but more complex situations could be tackled in a similar way.
Furthermore, we will focus on the mean field limit and will not discuss the models with a finite number of agents.

% subsec-subsec-subsec-subsec-subsec-subsec-subsec-subsec-subsec %
% subsec-subsec-subsec-subsec-subsec-subsec-subsec-subsec-subsec %
\subsection{Mathematical framework}
\label{sec:math-framework}

In the sequel, we will consider two populations with densities $(x,t) \mapsto m_{i}(x,t)$, $i=1,2$
and $x\in \RR^{n}.$ We note by $m$ the vector $(m_1, m_2)^*$, where the superscript $*$ denotes the transpose of a vector or a matrix. We note by $m_t$ the function $x \mapsto m(x,t)$. 
Feedback controls are functions $(x,m) \mapsto v_{i}(x,m)$ $\in \RR^{d}$ to be
chosen by the two populations. To simplify notation, we omit to write
explicitly a possible dependence on time. We note by $v$ the vector $(v_1, v_2)^*$.

For simplicity, we will assume that $m_{i}(\cdot,t) \in L^2(\RR^n)$. This allows us to use functional derivatives in the following sense. 
For a function $f : L^2(\RR^{n})^2 \to \RR^n$, we note by $\partial_{m_i} f$ the G{\^a}teaux derivative of $f$ with respect to the $i$-th density, so that for $m \in L^2(\RR^{n})^2, \tilde m \in L^2(\RR^{n})^2$,
\begin{equation}
\label{eq:def-grad-m}
	\frac{d}{d\theta} f(m_i+\theta\tilde m_i, m_{-i}) \big|_{\theta=0}
	= \int_{\R^n}\frac{\partial f}{\partial m_i}(m)(\xi) \, \tilde m_i(\xi)\,d\xi,
\end{equation}
where $i=1,2$ and $-i=2,1$ respectively.
This notion of differentiability will be sufficient for the purpose of this work. For a more rigorous treatment beyond the $L^2$ setting, one could rely on a more general notion of differentiability introduced by P.-L. Lions in his lectures at Coll{\`e}ge de France~\cite{PLL-CDF} and called $L$-derivative by Carmona and Delarue (see Chapter 5 in~\cite{MR3752669}).

\bigskip
The paper is organized as follows.
In Section \ref{sec:CMFC}, we consider the common mean field control problem with a single objective function. In some sense it is the simplest model to present because it does not involve any fixed point argument. By using calculus of variations we derive formally the adjoint equations and obtain a system of forward-backward PDEs characterizing the optimal solution. In Section~\ref{sec:NMFC}, we consider two McKean-Vlasov populations with their own objective functions. The Nash equilibrium cannot be described with PDEs in finite dimension and in this case one \emph{needs} to rely on a system of master equations.  
In Section~\ref{sec:MFG}, we study the corresponding mean field game settings, and obtain the adjoint equations for common or separated objective functions, respectively. 
In Section~\ref{sec:examples}, we provide two types of examples: we first revisit an example of crowd dynamics from~\cite{LachapelleWolfram-2011-MFG-congestion-aversion,MR3763083} and we then turn our attention to linear-quadratic models, for which we obtain systems of ODEs that can be related to Riccati equations.

% secsecsecsecsecsecsecsecsecsecsecsecsecsecsecsecsecsecsecsecsecsecsecsec %
% secsecsecsecsecsecsecsecsecsecsecsecsecsecsecsecsecsecsecsecsecsecsecsec %
% secsecsecsecsecsecsecsecsecsecsecsecsecsecsecsecsecsecsecsecsecsecsecsec %

\section{Common Mean Field Type Control}%{Two population Mean Field Control}
\label{sec:CMFC}
\setcounter{equation}{0}
\setcounter{theorem}{0}

% subsec-subsec-subsec-subsec-subsec-subsec-subsec-subsec-subsec %
% subsec-subsec-subsec-subsec-subsec-subsec-subsec-subsec-subsec %

\subsection{Definition of the problem} 
In this section, we investigate the situation where a global planner seeks to control in a distributed fashion two interacting populations driven by McKean-Vlasov dynamics, and tries to minimize a global cost. This setting can also be construed as a kind of social optimum for two populations (or two players with MKV dynamics) cooperating in order to minimize a cost which aggregates their objectives.
We call this problem \textbf{common mean field type control} (CMFC for short) and define it as follows.
 
 \begin{problem}[CMFC]\label{pb:CMFC}
 Find a feedback control $\hat v = (\hat v_1, \hat v_2)^*$ 
 minimizing the functional
\begin{align}
	J^{CMFC}(v_1,v_2)
	&=
	\sum_{i=1}^2  \left[ \int_0^T \int_{\RR^n} f_i(x,m^{v}_t,v_i(x,m^{v}_t)) m^{v}_i(x,t) dx \,dt + \int_{\RR^n} h_i(x,m^{v}_T) m^{v}_i(x,T) dx \right],
	\label{defJ-CMFC}
\end{align}
where $m^{v}= (m^{v}_1,m^{v}_2)^*$ solves the following system of Fokker-Planck (FP) equations: for $i=1,2,$
\begin{equation}\label{FPCeqGen}
	\frac{\partial m_i}{\partial t}(x,t)+A_i^*m_i(x,t)+\diver_x\, (g_i(x,m_t,v_i(x, m_t))m_i(x,t))) 
	=0 \, ,
	\qquad (x,t) \in \RR^n \times \RR_+ \, ,
\end{equation}
with initial condition $m_i(x,0)=\rho_{i,0}(x)$, $x \in \RR^n$.
\end{problem}
The functions
\begin{align}
	&g_i:\RR^n\times L^2(\RR^{n})^2\times\RR^{d}\rightarrow \RR^n \, , \quad (x,m,v_i) \mapsto g_i(x,m,v_i) \, ,
	\label{eq:CMFC-def-gi}
	\\
	&f_i: \RR^{n} \times L^2(\RR^{n})^2 \times\RR^{d} \rightarrow \RR \, , \quad (x,m,v_i) \mapsto f_i(x,m,v_i) \, ,
	\label{eq:CMFC-def-fi}
	\\
	&h_i: \RR^{n} \times L^2(\RR^{n})^2 \rightarrow \RR \, , \quad (x,m) \mapsto h_i(x,m) \, ,
	\label{eq:CMFC-def-hi}
\end{align}
are assumed to be differentiable with respect to all independent variables.
These functions as well as the control $v$ may also depend on time but we omit it to alleviate the notations. By $A_{i}^{*}$ we denote the formal dual operator of the differential operator 
\begin{equation}
	A_{i}\varphi(x)
	=
	-\sum_{\alpha,\beta=1}^{n}a_{i}^{\alpha\beta}(x)\dfrac{\partial^{2}\varphi(x)}{\partial x_{\alpha}\partial x_{\beta}} \, .
	\notag
\end{equation}
We assume sufficient smoothness on the drift functions $g_i$ as well as on
the feedback, to perform differentiation as needed. We note that in~\eqref{FPCeqGen} the
coupling holds only through the vector $m$.

One could consider more general cost functions (which can not be decomposed as a sum) but here we restrict our attention to the form given in~\eqref{defJ-CMFC} in order to facilitate the comparison with a setting where each population has its own cost (see Section~\ref{sec:NMFC}).

Although we are going to focus on PDE formulations in the rest of this work, let us mention that this model can also be motivated from a stochastic viewpoint as the optimal control of two McKean-Vlasov dynamics.
Indeed, one can consider a stochastic process $X^v=(X^{v}_1,X^{v}_2)$ in $\RR^{2n}$ 
with the following McKean-Vlasov dynamics
\begin{equation*}
	dX^{v}_i(t)=g_i\Big( \, X^{v}_i(t), \big(\mathcal L(X^{v}_1(t)), \mathcal L(X^{v}_2(t))\big)^*, v_i\left(X^{v}_i(t),\big(\mathcal L(X^{v}_1(t)), \mathcal L(X^{v}_2(t))\big)^*\Big) \,\right)\,dt
			+\sigma_i(X^{v}_i(t))\,dW_i(t)
\end{equation*}
for $i=1,2$, with $X^{v}(0)$ such that $X^{v}_i(0)$ has distribution $\rho_{i,0}$. Here, $W = (W_1, W_2)$ is a pair of independent $\RR^n$-valued Brownian motions on a probability space $(\Omega,\mathcal A, \mathcal P)$. Moreover $\mathcal L(X^{v}_i(t))$ denotes the distribution of $X^{v}_i(t)$. If we assume that these distributions have densities with respect to Lebesgue measure which are in $L^2(\RR^n)$, then these densities satisfy, at least formally, the FP equations~\eqref{FPCeqGen} with 
$
	a_i(x) = \frac12\sigma_i(x)\sigma_i^*(x).
$
Moreover, the objective functional~\eqref{defJ-CMFC} can be written as a sum of expectations, which correspond to the expected cost of each player.

\subsection{Necessary conditions of optimality}

We shall assume the existence of optimal feedbacks $(x,m) \mapsto \hat{v}_{i}(x,m)$
and look for necessary conditions of optimality. 
We denote by $m^{\hat v}=(m^{\hat v}_1,m^{\hat v}_2)^*$ the solutions of~\eqref{FPCeqGen} controlled by $\hat{v}=(\hat{v}_1,\hat{v}_2)^*$. 
We consider feedbacks $(x,m) \mapsto \hat{v}_{i}(x,m)+\theta\tilde{v}_{i}(x,m)$, $\theta \in \RR$, and we call $m_{i,\theta}$ the corresponding solutions of the
FP equations (\ref{FPCeqGen}). Then $\dfrac{m_{i,\theta}(x,t)-m^{\hat v}_{i}(x,t)}{\theta}\:\rightarrow\tilde{m}_{i}(x,t)$
pointwise as $\theta\rightarrow0,$ solution of
\begin{align}
	&\dfrac{\partial\tilde{m}_{i}}{\partial t}(x,t)
	+A_{i}^{*}\tilde{m}_{i}(x,t)
	+\diver_{x}\left(g_{i}(x,m^{\hat v}_t,\hat{v}_{i}(x,m^{\hat v}_t))\tilde{m}_{i}(x,t)  
	+ \left[
	\sum_{j=1}^{2}\int_{\RR^n}
			\left(\partial_{m_{j}}g_{i}(x,m^{\hat v}_t,\hat{v}_{i}(x,m^{\hat v}_t))(\xi)
			\vphantom{\dfrac{\partial g_{i}(x,m^{\hat v}_t,\hat{v}_{i}(x,m^{\hat v}_t))}{\partial v_{i}}}
			\right.\right.
			\right.
	\notag
	\\
	&
		\left.
		\left.\left.
				+\dfrac{\partial g_{i}(x,m^{\hat v}_t,\hat{v}_{i}(x,m^{\hat v}_t))}{\partial v_{i}}\partial_{m_{j}}\hat{v}_{i}(x,m^{\hat v}_t)(\xi)\right)
				\tilde{m}_{j}(\xi,t)d\xi
	\vphantom{\sum_{j=1}^{2}}
	+ \dfrac{\partial g_{i}(x,m^{\hat v}_t,\hat{v}_{i}(x,m^{\hat v}_t))}{\partial v_{i}}\tilde{v}_{i}(x,m^{\hat v}_t) \right] m^{\hat v}_{i}(x,t)\right)
	=0 \, ,
	\label{eq:CMFC-tildemi}
\end{align}
with the initial condition $\tilde{m}_{i}(x,0)=0$. Recall that $\partial_{m_{i}}\varphi(m)$ denotes the functional
derivative of a functional $(m_1,m_2)^* = m \mapsto \varphi(m)$ with respect to $m_i$, as defined by~\eqref{eq:def-grad-m}.  Note that $\tilde m$ depends on $\hat v$ and $ m^{\hat v}$, but we omit it to save notation.
We next
compute the G{\^a}teaux derivative of the functional $J^{CMFC}$ as follows
\begin{align}
	&\dfrac{d}{d\theta}J^{CMFC}(\hat{v}+\theta\tilde{v})|_{\theta=0}
	\label{eq:1-5}
	\\
	&=
	\sum_{i=1}^{2}\int_{0}^{T}\int_{\RR^n}f_{i}(x,m^{\hat v}_t,\hat{v}_{i}(x,m^{\hat v}_t))\tilde{m}_{i}(x,t)dxdt
	\notag
	\\
	&\quad+\sum_{i,j=1}^{2}\int_{0}^{T}\int_{\RR^n}\left[\int_{\RR^n}\left(\partial_{m_{j}}f_{i}(x,m^{\hat v}_t,\hat{v}_{i}(x,m^{\hat v}_t))(\xi)
	\vphantom{\dfrac{\partial f_{i}}{\partial v_{i}}}
	\right.\right.
	\notag
	\\
	& \qquad\qquad\qquad
	\left.\left.
		+\dfrac{\partial f_{i}(x,m^{\hat v}_t,\hat{v}_{i}(x,m^{\hat v}_t))}{\partial v_{i}}\partial_{m_{j}}\hat{v}_{i}(x,m^{\hat v}_t)(\xi)\right)\tilde{m}_{j}(\xi,t)d\xi\right]m^{\hat v}_{i}(x,t)dxdt
	\notag
	\\
	&\quad+\sum_{i=1}^{2}\int_{0}^{T}\int_{\RR^n}\dfrac{\partial f_{i}(x,m^{\hat v}_t,\hat{v}_{i}(x,m^{\hat v}_t))}{\partial v_{i}}\tilde{v}_{i}(x,m^{\hat v}_t)m^{\hat v}_{i}(x,t)dxdt
	\notag
	\\
	&\quad+\sum_{i=1}^{2}\int_{\RR^n}h_{i}(x,m^{\hat v}_{T})\tilde{m}_{i}(x,T)dx+\sum_{i,j=1}^{2}\int_{\RR^n}\left[\int_{\RR^n}\partial_{m_{j}}h_{i}(x,m^{\hat v}_{T})(\xi)\tilde{m}_{j}(\xi,T)d\xi\right]m^{\hat v}_{i}(x,T)dx \, .
	\notag
\end{align}
 We then introduce the functions $(x,t) \mapsto u^{\hat v}_{i}(x,t), i=1,2,$ solutions of the backward
equations
\begin{align}
	&-\dfrac{\partial u_{i}}{\partial t}(x,t)+A_{i}u_{i}(x,t)
	\label{eq:CMFC-HJBFPsys-HJB-TMP}
	\\
	= \, 
	&f_{i}(x,m^{\hat v}_t,\hat{v}_{i}(x,m^{\hat v}_t))+Du_{i}(x,t).g_{i}(x,m^{\hat v}_t,\hat{v}_{i}(x,m^{\hat v}_t))
	\notag
	\\
	&+\sum_{j=1}^{2}\int_{\RR^n}\left[
		\partial_{m_{i}}f_{j}(\xi,m^{\hat v}_t,\hat{v}_{j}(\xi,m^{\hat v}_t))(x)+\dfrac{\partial f_{j}(\xi,m^{\hat v}_t,\hat{v}_{j}(\xi,m^{\hat v}_t))}{\partial v_{j}}\partial_{m_{i}}\hat{v}_{j}(\xi,m^{\hat v}_t)(x)
		\right] m^{\hat v}_{j}(\xi,t)d\xi
	\notag
	\\
	&+\sum_{j=1}^{2}\int_{\RR^n}Du_{j}(\xi,t).\left(
			\partial_{m_{i}}g_{j}(\xi,m^{\hat v}_t,\hat{v}_{i}(\xi,m^{\hat v}_t))(x)+\dfrac{\partial g_{j}(\xi,m^{\hat v}_t,\hat{v}_{j}(\xi,m^{\hat v}_t))}{\partial v_{j}}\partial_{m_{i}}\hat{v}_{j}(\xi,m^{\hat v}_t)(x) 
			\right)m^{\hat v}_{j}(\xi,t)d\xi
	\notag
\end{align}
with terminal condition
$$
	u_{i}(x,T)
	=
	 h_{i}(x,m^{\hat v}_{T})+\sum_{j=1}^{2}\int_{\RR^n}\partial_{m_{i}}h_{j}(\xi,m^{\hat v}_{T})(x)m^{\hat v}_{j}(\xi,T) \, .
$$
Using (\ref{eq:CMFC-HJBFPsys-HJB-TMP}) in (\ref{eq:1-5}) it follows 
\begin{align*}
	&\dfrac{d}{d\theta}J^{CMFC}(\hat{v}+\theta\tilde{v})|_{\theta=0}
	\\
	=
	\, &\sum_{i=1}^{2}\int_{0}^{T}\int_{\RR^n}\tilde{m}_{i}(x,t)\left[-\dfrac{\partial u^{\hat v}_{i}}{\partial t}(x,t)+A_{i}u^{\hat v}_{i}(x,t)-Du^{\hat v}_{i}(x,t).g_{i}(x,m^{\hat v}_t,\hat{v}_{i}(x,m^{\hat v}_t))\vphantom{\sum_{j=1}^{2}}\right.
	\\
	&-\left.\sum_{j=1}^{2}\int_{\RR^n}Du^{\hat v}_{j}(\xi,t).\left(\partial_{m_{i}}g_{j}(\xi,m^{\hat v}_t,\hat{v}_{i}(\xi,m^{\hat v}_t))(x)
	\right.\right.
	\\
	&\qquad\qquad\qquad\qquad\qquad \left.\left.
	+\dfrac{\partial g_{j}(\xi,m^{\hat v}_t,\hat{v}_{j}(\xi,m^{\hat v}_t))}{\partial v_{j}}\partial_{m_{i}}\hat{v}_{j}(\xi,m^{\hat v}_t)(x)\right) m^{\hat v}_{j}(\xi,t)d\xi\right] dx dt
	\\
	&+ \sum_{i=1}^{2}\int_{0}^{T}\int_{\RR^n}\dfrac{\partial f_{i}(x,m^{\hat v}_t,\hat{v}_{i}(x,m^{\hat v}_t))}{\partial v_{i}}\tilde{v}_{i}(x,m^{\hat v}_t)m^{\hat v}_{i}(x,t)dxdt
		+\sum_{i=1}^{2}\int_{\RR^n}u^{\hat v}_{i}(x,T)\tilde{m}_{i}(x,T)dx \, .
\end{align*}
 Integrating by parts and using (\ref{eq:CMFC-tildemi}) after rearrangements,
we obtain
\begin{align}
	&\dfrac{d}{d\theta}J^{CMFC}(\hat{v}+\theta\tilde{v})|_{\theta=0}
	\notag
	\\
	= \, 
	&\sum_{i=1}^{2}\int_{0}^{T}\int_{\RR^n}
	\left[
		\dfrac{\partial f_{i}(x,m^{\hat v}_{t},\hat{v}_{i}(x,m^{\hat v}_{t}))}{\partial v_{i}}+Du^{\hat v}_{i}(x,t).\dfrac{\partial g_{i}(x,m^{\hat v}_{t},\hat{v}_{i}(x,m^{\hat v}_{t}))}{\partial v_{i}}
		\right]\tilde{v}_{i}(x,m^{\hat v}_{t})m^{\hat v}_{i}(x,t)dxdt \, .
	\notag
\end{align}
Since $\tilde{v}_{i}(x,m^{\hat v}_{t})$ can be an arbitrary function of $(x,t)$ and $m^{\hat v}_{i}(x,t)>0$, assuming the matrix $a_{i}^{\alpha\beta}(x)$
uniformly positive definite, we necessarily have 
\begin{equation}
	\dfrac{\partial f_{i}(x,m^{\hat v}_{t},\hat{v}_{i}(x,m^{\hat v}_{t}))}{\partial v_{i}}
	+Du^{\hat v}_{i}(x,t).\dfrac{\partial g_{i}(x,m^{\hat v}_{t},\hat{v}_{i}(x,m^{\hat v}_{t}))}{\partial v_{i}}
	=0,\,\text{a.e. }x,t \, .
	\label{eq:1-7}
\end{equation}
Let us introduce, for $i=1,2,$ the Hamiltonians 
\begin{equation}
	H_{i}(x,m,q_{i})
	=
	\inf_{v_{i}} \big[ f_{i}(x,m,v_{i})+q_{i}.g_{i}(x,m,v_{i}) \big] \, ,
	\label{eq:defH-CMFC}
\end{equation}
and the functions $(x,m,q_{i}) \mapsto \hat{\bfv}_{i}(x,m,q_{i})$ which achieve the infima, 
 that is 
$$
	H_{i}(x,m,q_{i}) = f_{i}\left(x,m,\hat{\bfv}_{i}(x,m,q_{i})\right)+q_{i}.g_{i}\left(x,m,\hat{\bfv}_{i}(x,m,q_{i})\right).
$$
We see from (\ref{eq:1-7}) that 
\begin{equation}
	\hat{v}_{i}(x,m^{\hat v}_{t}) = \hat{\bfv}_{i}(x,m^{\hat v}_{t},Du^{\hat v}_{i}(x,t)) \, .
		\label{eq:1-9}
\end{equation}
Hence we have
\begin{align}
	H_{i}(x,m^{\hat v}_{t},Du^{\hat v}_{i}(x,t))
	&=
	f_{i}(x,m^{\hat v}_{t},\hat{v}_{i}(x,m^{\hat v}_{t}))+Du^{\hat v}_{i}(x,t).g_{i}(x,m^{\hat v}_{t},\hat{v}_{i}(x,m^{\hat v}_{t})) \, ,
	\label{eq:1-10}
\end{align}
and, by~\eqref{eq:1-7} and~\eqref{eq:1-9},
\begin{align}
	\dfrac{\partial H_{i}}{\partial q_{i}}(x,m^{\hat v}_{t},Du^{\hat v}_{i}(x,t))
	&=
	g_{i}(x,m^{\hat v}_{t},\hat{v}_{i}(x,m^{\hat v}_{t})) \, .
	\label{eq:1-11}
\end{align}
Moreover,
\begin{align}
	&\sum_{j=1}^{2}\int_{\RR^n}\partial_{m_{i}}H_{j}(\xi,m^{\hat v}_{t},Du^{\hat v}_{j}(\xi,t))(x)m^{\hat v}_{j}(\xi,t)d\xi
	\label{eq:1-12}
	\\
	= \, &
	\sum_{j=1}^{2}\int_{\RR^n}\left[
				\partial_{m_{i}}f_{j}(\xi,m^{\hat v}_{t},\hat{v}_{j}(\xi,m^{\hat v}_{t}))(x)+\dfrac{\partial f_{j}(\xi,m^{\hat v}_{t},\hat{v}_{j}(\xi,m^{\hat v}_{t}))}{\partial v_{j}}\partial_{m_{i}}\hat{v}_{j}(\xi,m^{\hat v}_{t})(x)
				\right] m^{\hat v}_{j}(\xi,t)d\xi
	\notag
	\\
	& +\sum_{j=1}^{2}\int_{\RR^n}Du^{\hat v}_{j}(\xi,t).\left[
			\partial_{m_{i}}g_{j}(\xi,m^{\hat v}_{t},\hat{v}_{i}(\xi,m^{\hat v}_{t}))(x)+\dfrac{\partial g_{j}(\xi,m^{\hat v}_{t},\hat{v}_{j}(\xi,m^{\hat v}_{t}))}{\partial v_{j}}\partial_{m_{i}}\hat{v}_{j}(\xi,m^{\hat v}_{t})(x) \right] m^{\hat v}_{j}(\xi,t)d\xi \, .
	\notag
\end{align}
 Therefore $(u^{\hat v}_i, m^{\hat v}_i)_{i=1,2}$ solve the following system of HJB-FP equations
\begin{align}
	&-\dfrac{\partial u_{i}}{\partial t}(x,t) + A_{i}u_{i}(x,t)
	=
	H_{i}(x,m_{t},Du_{i}(x,t))+\sum_{j=1}^{2}\int_{\RR^n}\partial_{m_{i}}H_{j}(\xi,m_{t},Du_{j}(\xi,t))(x)m_{j}(\xi,t)d\xi
	\label{eq:CMFC-HJBFPsys-HJB}
	\\
	&\dfrac{\partial m_{i}}{\partial t}(x,t) 
	+ A_{i}^{*}m_{i}(x,t) 
	+ \diver_{x}\left(\dfrac{\partial H_{i}}{\partial q_{i}}(x,m_{t},Du_{i}(x,t))m_{i}(x,t)\right)=0
	\label{eq:CMFC-HJBFPsys-FP}
\end{align}
with terminal and initial condition
$$
	u_{i}(x,T)=h_{i}(x,m_{T})+\sum_{j=1}^{2}\int_{\RR^n}\partial_{m_{i}}h_{j}(\xi,m_{T})(x)m_{j}(\xi,T)d\xi \, ,
	\qquad m_{i}(x,0)=\rho_{i,0}(x) \, ,
$$
where $m = (m_1,m_2)^*$ denotes the vector of solutions of~\eqref{eq:CMFC-HJBFPsys-FP}.

The PDE system~\eqref{eq:CMFC-HJBFPsys-HJB}--\eqref{eq:CMFC-HJBFPsys-FP} extends to two populations the PDE system for mean field control with a single population (see e.g.~\cite{MR3134900}, Chapter 4, page 18).  In Section~\ref{sec:NMFC}, we will present a different way to extend to two populations the mean field control framework.

\begin{remark}
It is important to notice that finding the solution $(u_{i},m_{i})_{i=1,2}$ of the above PDE system allows to compute the functions $(x,t) \mapsto \hat{v}_{i}(x,m_{t})$
but not the feedbacks $(x,m) \mapsto \hat{v}_{i}(x,m)$. In other words, the optimal controls can be computed only along the optimal flows of distributions (that is, the solution $t \mapsto m_t$ obtained by solving the PDE system), but $\hat v_i(x,m)$ is not known for all possible $m$. This can be known only through the master equations, as we explain below. 
\end{remark}

\subsection{Master and Bellman equations}

The notion of master equation (for a single population) has been introduced by P.-L. Lions in the context of mean field games~\cite{PLL-CDF} and has been studied e.g. in~\cite{MR3343705,MR3652408,MR3332710}. For more details, the reader is referred to~\cite{MR3753660}.

This section is devoted to the formal introduction of an analogous equation (or, rather, a system of analogous equations) for two-population CMFC. 
In this setting, the master equations are equations for functions $(x,m,t) \mapsto U_{i}(x,m,t)$ such
that, in particular, for $m$ solving~\eqref{eq:CMFC-HJBFPsys-FP}, there holds 
\begin{equation}\label{eq:rel-ui-Ui}
	u_{i}(x,t)=U_{i}(x,m_{t},t).
\end{equation}
These equations are self-contained, whereas
the HJB equations~\eqref{eq:CMFC-HJBFPsys-HJB} for $u_{i}$ are not since they need to be coupled
with the functions $m_{i}$ solutions of the FP equations~\eqref{eq:CMFC-HJBFPsys-FP}. 
Assuming~\eqref{eq:rel-ui-Ui} holds, using the FP equation~\eqref{eq:CMFC-HJBFPsys-FP} and integration by parts, we obtain
\begin{align*}
	\dfrac{\partial u_{i}}{\partial t}(x,t)
	&=
	\dfrac{\partial U_{i}}{\partial t}(x,m_{t},t)
	+\sum_{j=1}^{2}\int_{\RR^{n}}\partial_{m_{j}}U_{i}(x,m_{t},t)(\xi)\dfrac{\partial m_{j}}{\partial t}(\xi,t)d\xi
	\\
	&=
	\dfrac{\partial U_{i}}{\partial t}(x,m_{t},t)
	+\sum_{j=1}^{2}\int_{\RR^{n}} \left[ A_{j\xi} + \dfrac{\partial H_{j}}{\partial q_{j}}(x,m_{t},Du_{j}(x,t)) D_{\xi}\right]\partial_{m_{j}}U_{i}(x,m_{t},t)(\xi)  m_{j}(\xi,t) d\xi \, .
\end{align*}
Here, $\partial_{m_{j}}U_{i}(x,m,t)$ denotes a derivative in the sense of~\eqref{eq:def-grad-m}; it is a function of the space variable and $D_{\xi}\partial_{m_{j}}U_{i}(x,m,t)(\xi)$ denotes its gradient. The notation $A_{i\xi}$ is to be understood in a similar way.
 Using the adjoint equations~\eqref{eq:CMFC-HJBFPsys-HJB}, we identify the master equations 
\begin{align}
	&-\dfrac{\partial U_{i}}{\partial t}(x,m,t) + A_{ix}U_{i}(x,m,t)
	\label{eq:1-14}
	\\
	= \, 
	&-\sum_{j=1}^{2}\int_{\RR^{n}}A_{j\xi}\partial_{m_{j}}U_{i}(x,m,t)(\xi)\,m_{j}(\xi)d\xi
	\notag
	\\
	&+H_{i}(x,m,DU_{i}(x,m,t))+\sum_{j=1}^{2}\int_{\RR^{n}}\partial_{m_{i}}H_{j}(\xi,m_{t},DU_{j}(\xi,m,t))(x)m_{j}(\xi,t)d\xi
	\notag
	\\
	&+\sum_{j=1}^{2}\int_{\RR^{n}}D_{\xi}\partial_{m_{j}}U_{i}(x,m,t)(\xi)\dfrac{\partial H_{j}}{\partial q_{j}}(\xi,m,DU_{j}(\xi,m,t))m_{j}(\xi)d\xi \, ,
	\notag
\end{align}
with terminal condition
$$
	U_{i}(x,m,T)=h_{i}(x,m)+\sum_{j=1}^{2}\int_{\RR^{n}}\partial_{m_{i}}h_{j}(\xi,m)(x)m_{j}(\xi)d\xi \, .
$$

 We then identify $U_{i}$ as a functional derivative in $m_{i}$, namely 
\begin{equation}
	U_{i}(x,m,t)
	=\partial_{m_{i}} V(m,t)(x) \, ,
	\label{eq:1-15}
\end{equation}
with $(m,t)\mapsto V(m,t)$ solution of the Bellman equation
\begin{align}
	-\dfrac{\partial V}{\partial t}(m,t)
	+\sum_{j=1}^{2}\int_{\RR^{n}}A_{j}\partial_{m_{j}}V(m,t)(x)m_{j}(x)dx
	&= \sum_{j=1}^{2}\int_{\RR^{n}}H_{j}(x,m,\partial_{m_{j}}V(m,t)(x)) m_j(x)dx\, ,
	\label{eq:1-16}
\end{align}
with terminal condition
$$
	V(m,T)
	=\sum_{j=1}^{2}\int_{\RR^{n}}h_{j}(x,m)m_{j}(x)dx \, .
$$
 Thanks to the functions $U_{i}$ we can completely characterize
the feedbacks by (recall that $\hat{\bfv}_{i}$ minimizes the Hamiltonian defined by~\ref{eq:defH-CMFC})
\[
	\hat{v}_{i}(x,m)=\hat{\bfv}_{i}(x,m,DU_{i}(x,m,t)) \, .
\]
\begin{remark}
From (\ref{eq:1-15}) we can assert that 
\begin{equation}
	\partial_{m_{j}}U_{i}(x,m,t)(\xi)=\partial_{m_{i}}U_{j}(\xi,m,t)(x)
	=
	\partial_{m_{i} m_{j}}^{2}V(m,t)(x,\xi) \, .
	\label{eq:1-17}
\end{equation}
\end{remark}

\begin{remark}
	The Bellman equation~\eqref{eq:1-16} could also be obtained directly by a dynamic programming argument similarly to the case of single population, see e.g.~\cite{MR3332710,MR3343705,MR3501391}. 
\end{remark}

% secsecsecsecsecsecsecsecsecsecsecsecsecsecsecsecsecsecsecsecsecsecsecsec %
% secsecsecsecsecsecsecsecsecsecsecsecsecsecsecsecsecsecsecsecsecsecsecsec %
% secsecsecsecsecsecsecsecsecsecsecsecsecsecsecsecsecsecsecsecsecsecsecsec %

\section{Nash Mean Field Type Control Problem\label{sec:NMFC,-NASH-MEAN} }
\label{sec:NMFC}

\subsection{Definition of the problem}

In this section, we consider a situation where the agents among each population cooperate, but the two populations compete and we
try to find a Nash equilibrium between them. 
A key point is that, here again, each population chooses
a feedback $(x,m) \mapsto v_{i}(x,m)$, so we incorporate both $m_{1}$ and $m_{2}$ in
the feedback as in the case CMFC. The FP equations describing the
evolution of $m_{i}(x,t)$ are still given by~\eqref{FPCeqGen}. However
there is not a common cost functional. Instead, each population has its own functional. We call this problem \textbf{Nash mean field control} (NMFC for short) and define it as follows.

\begin{problem}[NMFC] Find a Nash equilibrium $\hat v = (\hat v_1, \hat v_2)^*$ 
for the cost functionals
\begin{equation}
	J^{NMFC}_{i}(v_{1},v_{2})
	=
	\int_{0}^{T}\int_{\RR^{n}}f_{i}(x,m^v_{t},v_{i}(x,m^v_{t})) m^v_{i}(x,t) dxdt+\int_{\RR^{n}}h_{i}(x,m^v_{T})m^v_{i}(x,T) dx\label{eq:2-1}
\end{equation}
where $m^{v}= (m^{v}_1,m^{v}_2)^*$ satisfies the PDEs~\eqref{FPCeqGen}
with initial conditions $m_i(\cdot,0)=\rho_{i,0}$.
\end{problem}
Here again, $g_i,$ $f_i$ and $h_i$ are as in~\eqref{eq:CMFC-def-gi}, \eqref{eq:CMFC-def-fi} and~\eqref{eq:CMFC-def-hi} respectively, hence $J^{CMFC} = \sum_{i=1}^2 J^{NMFC}_{i}.$

\begin{remark}
	This situation can also be viewed as a competition between two players, each having a dynamics of McKean-Vlasov type. 
\end{remark}

In the sequel $\hat{v} = (\hat{v}_1,\hat{v}_2)^*$ represents a solution to this problem (assuming it exists). In other words, we have for all $v = (v_{1},v_{2})^*$
$$
	J^{NMFC}_{1}(\hat v_{1}, \hat v_{2}) \leq J^{NMFC}_{1}(v_{1}, \hat v_{2}) 
	\quad \hbox{and} \quad
	J^{NMFC}_{2}(\hat v_{1}, \hat v_{2}) \leq J^{NMFC}_{2}(\hat v_{1},  v_{2}) \, .
$$
We shall write necessary conditions for the existence of such an equilibrium $\hat{v}$.

% subsec-subsec-subsec-subsec-subsec-subsec-subsec-subsec-subsec %
% subsec-subsec-subsec-subsec-subsec-subsec-subsec-subsec-subsec %
\subsection{Problem of player $1$}

For player $1,$ $\hat{v}_{2}$ is a fixed function and $\hat{v}_{1}$
solves the following (one-population) MFC problem.

\begin{problem}[NMFC: Problem of player $1$]
\label{pb:NMFC-player1}
 Minimize
\begin{align}
	&\cJ^{NMFC}_{1}(v_{1}) = J^{NMFC}_{1}(v_{1}, \hat v_{2})
	\notag
	\\
	=\,
	& \int_{0}^{T}\int_{\RR^{n}}f_{1}(x,m^{v_{1}, \hat v_{2}}_{t},v_{1}(x,m^{v_{1}, \hat v_{2}}_{t})) m^{v_{1}, \hat v_{2}}_1(x,t) dxdt+\int_{\RR^{n}}h_{1}(x,m^{v_{1}, \hat v_{2}}_{T}) m^{v_{1}, \hat v_{2}}_1(x,T) dx
	\label{eq:2-3}
\end{align}
where $m^{^{v_{1}, \hat v_{2}}} = \left(m^{^{v_{1}, \hat v_{2}}}_1, m^{^{v_{1}, \hat v_{2}}}_2\right)^*$ solves~\eqref{FPCeqGen} controlled by $(v_{1}, \hat v_{2})$, that is,
\begin{align}
	&\dfrac{\partial m_{1}}{\partial t}(x,t) + A_{1}^{*}m_{1}(x,t) + \diver_{x}(g_{1}(x,m_{t},v_{1}(x,m_{t}))m_{1}(x,t))=0
	\label{eq:NMFC-FP1-1}
	\\
	&\dfrac{\partial m_{2}}{\partial t}(x,t) + A_{2}^{*}m_{2}(x,t) + \diver_{x}(g_{2}(x,m_{t},\hat{v}_{2}(x,m_{t}))m_{2}(x,t))=0
	\label{eq:NMFC-FP1-2}
\end{align}
with the initial conditions $m_{i}(\cdot,0)=\rho_{i,0}$.
\end{problem}
In this section, we denote by $m^{\hat v}=(m^{\hat v}_1,m^{\hat v}_2)^*$ the solutions of~\eqref{eq:NMFC-FP1-1}--\eqref{eq:NMFC-FP1-2} when player $1$ chooses $\hat{v}_1$ as a control. 
Consider a feedback
	$(x,m) \mapsto \hat{v}_{1}(x,m)+\theta\tilde{v}_{1}(x,m)$ and denote by $m_{1,\theta}$,
$m_{2,\theta}$ the corresponding solutions of the FP equations~\eqref{eq:NMFC-FP1-1},~\eqref{eq:NMFC-FP1-2} respectively. Then 
\[
	\dfrac{m_{1,\theta}-m^{\hat v}_{1}}{\theta}(x,t) \xrightarrow[\theta \to 0]{} \tilde{m}_{1}(x,t) 
	\, , 
	\qquad \dfrac{m_{2,\theta}-m^{\hat v}_{2}}{\theta}(x,t) \xrightarrow[\theta \to 0]{} \tilde{m}_{-1}(x,t) \, .
\]
 Notice that we have written $\tilde{m}_{-1}(x,t)$ and not $\tilde{m}_{2}(x,t)$
for the second limit. Indeed, we are considering the problem of player~$1,$
and the second limit describes the impact of his choice on the second
FP equation~\eqref{eq:NMFC-FP1-2}. So the index~$1$ characterizes the first player, and
the sign~$\pm$ refers to his FP equation or the FP equation of
the opponent (player $2$). We obtain for $\tilde m_{1}$
\begin{align}
	&\dfrac{\partial\tilde{m}_{1}}{\partial t}(x,t)
	+A_{1}^{*}\tilde{m}_{1}(x,t)
	+\diver_{x} \left\{
	g_{1}(x,m^{\hat v}_{t},\hat{v}_{1}(x,m^{\hat v}_{t}))\tilde{m}_{1}(x,t) 
	\vphantom{\dfrac{\partial g_{1}(x,m^{\hat v}_{t},\hat{v}_{1}(x,m^{\hat v}_{t}))}{\partial v_{1}}} \right.
	\label{eq:NMFC-tildem1}
	\\
	&
	+ \left[\int_{\RR^{n}}\left(\partial_{m_{1}}g_{1}(x,m^{\hat v}_{t},\hat{v}_{1}(x,m^{\hat v}_{t}))(\xi)
	+\dfrac{\partial g_{1}(x,m^{\hat v}_{t},\hat{v}_{1}(x,m^{\hat v}_{t}))}{\partial v_{1}}\partial_{m_{1}}\hat{v}_{1}(x,m^{\hat v}_{t})(\xi)\right)\tilde{m}_{1}(\xi,t)d\xi
	\right.
	\notag
	\\
	&
	\quad + 
	\int_{\RR^{n}}\left(\partial_{m_{2}}g_{1}(x,m^{\hat v}_{t},\hat{v}_{1}(x,m^{\hat v}_{t}))(\xi)
	+\dfrac{\partial g_{1}(x,m^{\hat v}_{t},\hat{v}_{1}(x,m^{\hat v}_{t}))}{\partial v_{1}}\partial_{m_{2}}\hat{v}_{1}(x,m^{\hat v}_{t})(\xi)\right)\tilde{m}_{-1}(\xi,t)d\xi 
	\notag
	\\
	&
	\left. \left.
	\quad + \dfrac{\partial g_{1}(x,m^{\hat v}_{t},\hat{v}_{1}(x,m^{\hat v}_{t}))}{\partial v_{1}}\tilde{v}_{1}(x,m^{\hat v}_{t}) 
	\vphantom{\dfrac{\partial g_{1}(x,m^{\hat v}_{t},\hat{v}_{1}(x,m^{\hat v}_{t}))}{\partial v_{1}}} 
	\right] m^{\hat v}_{1}(x,t)
	\right \}
	= 0 \, ,
	\notag
\end{align}
 and for $\tilde m_{-1}$
\begin{align}
	&\dfrac{\partial\tilde{m}_{-1}}{\partial t}(x,t) + A_{2}^{*}\tilde{m}_{-1}(x,t)
	+\diver_{x} \left\{ 
		g_{2}(x,m^{\hat v}_{t},\hat{v}_{2}(x,m^{\hat v}_{t}))\tilde{m}_{-1}(x,t) %\right)
		\vphantom{\dfrac{\partial g_{2}(x,m^{\hat v}_{t},\hat{v}_{2}(x,m^{\hat v}_{t}))}{\partial v_{2}}}
		\right.
	\label{eq:NMFC-tildem-1}
	\\
	&%+\diver_{x}\left( 
	+ \left[\int_{\RR^{n}}\left(\partial_{m_{1}}g_{2}(x,m^{\hat v}_{t},\hat{v}_{2}(x,m^{\hat v}_{t}))(\xi)
	+\dfrac{\partial g_{2}(x,m^{\hat v}_{t},\hat{v}_{2}(x,m^{\hat v}_{t}))}{\partial v_{2}}\partial_{m_{1}}\hat{v}_{2}(x,m^{\hat v}_{t})(\xi)\right)\tilde{m}_{1}(\xi,t)d\xi %\right]m^{\hat v}_{2}(x,t)%\right)
	\right.
	\notag
	\\
	&%+\diver_{x}\left( 
%	\left.\left.
	\quad + %\left[
	\int_{\RR^{n}}\left(\partial_{m_{2}}g_{2}(x,m^{\hat v}_{t},\hat{v}_{2}(x,m^{\hat v}_{t}))(\xi)
	\vphantom{\dfrac{\partial g_{2}(x,m^{\hat v}_{t},\hat{v}_{2}(x,m^{\hat v}_{t}))}{\partial v_{2}}}
	\right.
	\notag
	\\
	&\qquad\qquad
	\left.\left.\left.
	+\dfrac{\partial g_{2}(x,m^{\hat v}_{t},\hat{v}_{2}(x,m^{\hat v}_{t}))}{\partial v_{2}}\partial_{m_{2}}\hat{v}_{2}(x,m^{\hat v}_{t})(\xi)\right)\tilde{m}_{-1}(\xi,t)d\xi
	\right] m^{\hat v}_{2}(x,t)
	\right\}
	= 0 \, ,
	\notag
\end{align}
with the initial conditions
$$
	\tilde{m}_{1}(x,0)=0, \qquad \tilde{m}_{-1}(x,0)=0 \, .
$$
Similarly to the CMFC case, we can compute 
\begin{align}
	 &\dfrac{d}{d\theta}\cJ^{NMFC}_{1}(\hat{v}_{1}+\theta\tilde{v}_{1})|_{\theta=0}
	 \label{eq:NMFC-FOC-1}
	 \\
	 &=\int_{0}^{T}\int_{\RR^{n}}f_{1}(x,m^{\hat v}_{t},\hat{v}_{1}(x,m^{\hat v}_{t}))\tilde{m}_{1}(x,t)dxdt
	 \notag
	 \\
	 & +\int_{0}^{T}\int_{\RR^{n}}\left[ 
	 	\int_{\RR^{n}}\left(\partial_{m_{1}}f_{1}(x,m^{\hat v}_{t},\hat{v}_{1}(x,m^{\hat v}_{t}))(\xi)
		\vphantom{\dfrac{\partial f_{1}(x,m^{\hat v}_{t},\hat{v}_{1}(x,m^{\hat v}_{t}))}{\partial v_{1}}}
		\right.\right.
	\notag
	\\
	&\qquad\qquad
	\left.\left.
	 	 +\dfrac{\partial f_{1}(x,m^{\hat v}_{t},\hat{v}_{1}(x,m^{\hat v}_{t}))}{\partial v_{1}}\partial_{m_{1}}\hat{v}_{1}(x,m^{\hat v}_{t})(\xi)\right)\tilde{m}_{1}(\xi,t)d\xi 
		 \right] m^{\hat v}_{1}(x,t)dxdt
	 \notag
	 \\
	 & +\int_{0}^{T}\int_{\RR^{n}}\left[
	 	\int_{\RR^{n}}\left(\partial_{m_{2}}f_{1}(x,m^{\hat v}_{t},\hat{v}_{1}(x,m^{\hat v}_{t}))(\xi)
		\vphantom{\dfrac{\partial f_{1}(x,m^{\hat v}_{t},\hat{v}_{1}(x,m^{\hat v}_{t}))}{\partial v_{1}}}
		\right.\right.
	\notag
	\\
	&\qquad\qquad
	\left.\left.
	 	+\dfrac{\partial f_{1}(x,m^{\hat v}_{t},\hat{v}_{1}(x,m^{\hat v}_{t}))}{\partial v_{1}}\partial_{m_{2}}\hat{v}_{1}(x,m^{\hat v}_{t})(\xi)\right)\tilde{m}_{-1}(\xi,t)d\xi
		\right] m^{\hat v}_{1}(x,t)dxdt
	 \notag
	 \\
	 & +\int_{0}^{T}\int_{\RR^{n}}\dfrac{\partial f_{1}(x,m^{\hat v}_{t},\hat{v}_{1}(x,m^{\hat v}_{t}))}{\partial v_{1}}\tilde{v}_{1}(x,m^{\hat v}_{t})m^{\hat v}_{1}(x,t)dxdt
	 \notag
	 \\
	 & +\int_{\RR^{n}}h_{1}(x,m^{\hat v}_{T})\tilde{m}_{1}(x,T)dx
	 \notag
	 \\
	 &+\int_{\RR^{n}}\left[
	 	\int_{\RR^{n}}\left(\partial_{m_{1}}h_{1}(x,m^{\hat v}_{T})(\xi)\tilde{m}_{1}(\xi,T)+\partial_{m_{2}}h_{1}(x,m^{\hat v}_{T})(\xi)\tilde{m}_{-1}(\xi,T)\right)d\xi
		\right] m^{\hat v}_{1}(x,T)dx \, .
	 \notag
\end{align}

% subsec-subsec-subsec-subsec-subsec-subsec-subsec-subsec-subsec %
% subsec-subsec-subsec-subsec-subsec-subsec-subsec-subsec-subsec %

\subsection{HJB-FP system for player $1$}

We introduce the functions $u^{\hat v}_{1},$ and $u^{\hat v}_{-1}$, solutions of
the equations 
\begin{align}
	& -\dfrac{\partial u_{1}}{\partial t}(x,t) + A_{1}u_{1}(x,t)
	\label{eq:NMFC-u1-hat}
	\\
	& =
	f_{1}(x,m^{\hat v}_{t},\hat{v}_{1}(x,m^{\hat v}_{t}))+Du_{1}(x,t).g_{1}(x,m^{\hat v}_{t},\hat{v}_{1}(x,m^{\hat v}_{t}))
	\notag
	\\
	&\qquad +\int_{\RR^{n}}\left[\partial_{m_{1}}f_{1}(\xi,m^{\hat v}_{t},\hat{v}_{1}(\xi,m^{\hat v}_{t}))(x)
		+\dfrac{\partial f_{1}(\xi,m^{\hat v}_{t},\hat{v}_{1}(\xi,m^{\hat v}_{t}))}{\partial v_{1}}\partial_{m_{1}}\hat{v}_{1}(\xi,m^{\hat v}_{t})(x)\right]m^{\hat v}_{1}(\xi,t)d\xi
	\notag
	\\
	&\qquad +\int_{\RR^{n}}\left[ Du_{1}(\xi,t).\left(
	\partial_{m_{1}}g_{1}(\xi,m^{\hat v}_{t},\hat{v}_{1}(\xi,m^{\hat v}_{t}))(x)
	\vphantom{\dfrac{\partial g_{1}(\xi,m^{\hat v}_{t},\hat{v}_{1}(\xi,m^{\hat v}_{t}))}{\partial v_{1}}}
	\right.\right.
	\notag
	\\
	&\qquad\qquad\qquad\qquad\qquad
	\left.\left.
		+\dfrac{\partial g_{1}(\xi,m^{\hat v}_{t},\hat{v}_{1}(\xi,m^{\hat v}_{t}))}{\partial v_{1}}\partial_{m_{1}}\hat{v}_{1}(\xi,m^{\hat v}_{t})(x)
		\right) \right]m^{\hat v}_{1}(\xi,t)d\xi
	\notag
	\\
	&\qquad +\int_{\RR^{n}}\left[ Du_{-1}(\xi,t).\left(
		\partial_{m_{1}}g_{2}(\xi,m^{\hat v}_{t},\hat{v}_{2}(\xi,m^{\hat v}_{t}))(x)
	\vphantom{\dfrac{\partial g_{1}(\xi,m^{\hat v}_{t},\hat{v}_{1}(\xi,m^{\hat v}_{t}))}{\partial v_{1}}}
	\right.\right.
	\notag
	\\
	&\qquad\qquad\qquad\qquad\qquad
	\left.\left.
		+\dfrac{\partial g_{2}(\xi,m^{\hat v}_{t},\hat{v}_{2}(\xi,m^{\hat v}_{t}))}{\partial v_{2}}\partial_{m_{1}}\hat{v}_{2}(\xi,m^{\hat v}_{t})(x)
		\right)\right] m^{\hat v}_{2}(\xi,t)d\xi
	\notag
\end{align}
with terminal condition
$$
	u_{1}(x,T)=h_{1}(x,m^{\hat v}_{T})+\int_{\RR^{n}}\partial_{m_{1}}h_{1}(\xi,m^{\hat v}_{T})(x)m^{\hat v}_{1}(\xi,T)d\xi \, ,
$$
and 
\begin{align}
	&-\dfrac{\partial u_{-1}}{\partial t}(x,t) + A_{2}u_{-1}(x,t)
	\label{eq:NMFC-u-1-hat}
	\\
	& = Du_{-1}(x,t).g_{2}(x,m^{\hat v}_{t},\hat{v}_{2}(x,m^{\hat v}_{t}))
	\notag
	\\
	& \qquad +\int_{\RR^{n}} \left[ \partial_{m_{2}}f_{1}(\xi,m^{\hat v}_{t},\hat{v}_{1}(\xi,m^{\hat v}_{t}))(x)
		+\dfrac{\partial f_{1}(\xi,m^{\hat v}_{t},\hat{v}_{1}(\xi,m^{\hat v}_{t}))}{\partial v_{1}}\partial_{m_{2}}\hat{v}_{1}(\xi,m^{\hat v}_{t})(x)\right] m^{\hat v}_{1}(\xi,t)d\xi
	\notag
	\\
	& \qquad +\int_{\RR^{n}} \left[ Du_{1}(\xi,t).\left(
		\partial_{m_{2}}g_{1}(\xi,m^{\hat v}_{t},\hat{v}_{1}(\xi,m^{\hat v}_{t}))(x)
	\vphantom{\dfrac{\partial g_{1}(\xi,m^{\hat v}_{t},\hat{v}_{1}(\xi,m^{\hat v}_{t}))}{\partial v_{1}}}
	\right.\right.
	\notag
	\\
	&\qquad\qquad\qquad\qquad\qquad
	\left.\left.
		+\dfrac{\partial g_{1}(\xi,m^{\hat v}_{t},\hat{v}_{1}(\xi,m^{\hat v}_{t}))}{\partial v_{1}}\partial_{m_{2}}\hat{v}_{1}(\xi,m^{\hat v}_{t})(x)
		\right) \right] m^{\hat v}_{1}(\xi,t)d\xi
	\notag
	\\
	& \qquad +\int_{\RR^{n}} \left[ Du_{-1}(\xi,t).\left(
		\partial_{m_{2}}g_{2}(\xi,m^{\hat v}_{t},\hat{v}_{2}(\xi,m^{\hat v}_{t}))(x)
	\vphantom{\dfrac{\partial g_{1}(\xi,m^{\hat v}_{t},\hat{v}_{1}(\xi,m^{\hat v}_{t}))}{\partial v_{1}}}
	\right.\right.
	\notag
	\\
	&\qquad\qquad\qquad\qquad\qquad
	\left.\left.
		+\dfrac{\partial g_{2}(\xi,m^{\hat v}_{t},\hat{v}_{2}(\xi,m^{\hat v}_{t}))}{\partial v_{2}}\partial_{m_{2}}\hat{v}_{2}(\xi,m^{\hat v}_{t})(x)
		\right) \right] m^{\hat v}_{2}(\xi,t)d\xi
	\notag
\end{align}
with terminal condition
$$
	u_{-1}(x,T)=\int_{\RR^{n}}\partial_{m_{2}}h_{1}(\xi,m^{\hat v}_{T})(x)m^{\hat v}_{1}(\xi,T)d\xi \, .
$$
 We now use (\ref{eq:NMFC-u1-hat}), (\ref{eq:NMFC-u-1-hat}) in (\ref{eq:NMFC-FOC-1}) to
obtain 
\begin{align*}
	&\dfrac{d}{d\theta}\cJ^{NMFC}_{1}(\hat{v}_{1}+\theta\tilde{v}_{1})|_{\theta=0}
	\\
	=
	\, &\int_{0}^{T}\int_{\RR^{n}}\tilde{m}_{1}(x,t)
	\left[
		-\dfrac{\partial u^{\hat v}_{1}}{\partial t}(x,t) + A_{1}u^{\hat v}_{1}(x,t) - Du^{\hat v}_{1}(x,t).g_{1}(x,m^{\hat v}_{t},\hat{v}_{1}(x,m^{\hat v}_{t}))
	\right.
	\\
	&\quad 
		-\int_{\RR^{n}}\left\{ Du^{\hat v}_{1}(\xi,t).\left(
			\partial_{m_{1}}g_{1}(\xi,m^{\hat v}_{t},\hat{v}_{1}(\xi,m^{\hat v}_{t}))(x)
			\vphantom{\dfrac{\partial g_{1}(\xi,m^{\hat v}_{t},\hat{v}_{1}(\xi,m^{\hat v}_{t}))}{\partial v_{1}}}
		\right.\right.
	\notag
	\\
	&\qquad\qquad \left.\left.
			+\dfrac{\partial g_{1}(\xi,m^{\hat v}_{t},\hat{v}_{1}(\xi,m^{\hat v}_{t}))}{\partial v_{1}}\partial_{m_{1}}\hat{v}_{1}(\xi,m^{\hat v}_{t})(x)
		\right) \right\} m^{\hat v}_{1}(\xi,t)d\xi
	\\
	&\quad 
	\left.
		-\int_{\RR^{n}} \left\{ Du^{\hat v}_{-1}(\xi,t).\left(
			\partial_{m_{1}}g_{2}(\xi,m^{\hat v}_{t},\hat{v}_{2}(\xi,m^{\hat v}_{t}))(x)
			\vphantom{\dfrac{\partial g_{1}(\xi,m^{\hat v}_{t},\hat{v}_{1}(\xi,m^{\hat v}_{t}))}{\partial v_{1}}}
		\right.\right.\right.
	\notag
	\\
	& \qquad\qquad \left.\left.\left.
			+ \dfrac{\partial g_{2}(\xi,m^{\hat v}_{t},\hat{v}_{2}(\xi,m^{\hat v}_{t}))}{\partial v_{2}}\partial_{m_{1}}\hat{v}_{2}(\xi,m^{\hat v}_{t})(x)
			\right) \right\} m^{\hat v}_{2}(\xi,t)d\xi
	\right]dxdt +
	\allowdisplaybreaks
	\\
	& 
	+\int_{0}^{T}\int_{\RR^{n}}\tilde{m}_{-1}(x,t)
	\left[
		-\dfrac{\partial u^{\hat v}_{-1}}{\partial t}(x,t) + A_{2}u^{\hat v}_{-1}(x,t) - Du^{\hat v}_{-1}(x,t).g_{2}(x,m^{\hat v}_{t},\hat{v}_{2}(x,m^{\hat v}_{t}))
	\right.
	\\
	& \quad
		-\int_{\RR^{n}} \left\{ Du^{\hat v}_{1}(\xi,t).\left(
			\partial_{m_{2}}g_{1}(\xi,m^{\hat v}_{t},\hat{v}_{1}(\xi,m^{\hat v}_{t}))(x)
			\vphantom{\dfrac{\partial g_{1}(\xi,m^{\hat v}_{t},\hat{v}_{1}(\xi,m^{\hat v}_{t}))}{\partial v_{1}}}
		\right.\right.
	\notag
	\\
	& \qquad\qquad \left.\left.
			+ \dfrac{\partial g_{1}(\xi,m^{\hat v}_{t},\hat{v}_{1}(\xi,m^{\hat v}_{t}))}{\partial v_{1}}\partial_{m_{2}}\hat{v}_{1}(\xi,m^{\hat v}_{t})(x)
			\right)
		\right\} m^{\hat v}_{1}(\xi,t)d\xi
	\\
	&\quad 
	\left.
		-\int_{\RR^{n}} \left\{ Du^{\hat v}_{-1}(\xi,t).\left(
			\partial_{m_{2}}g_{2}(\xi,m^{\hat v}_{t},\hat{v}_{2}(\xi,m^{\hat v}_{t}))(x)
			\vphantom{\dfrac{\partial g_{1}(\xi,m^{\hat v}_{t},\hat{v}_{1}(\xi,m^{\hat v}_{t}))}{\partial v_{1}}}
		\right.\right.\right.
	\notag
	\\
	&\qquad\qquad \left.\left.\left.
			+\dfrac{\partial g_{2}(\xi,m^{\hat v}_{t},\hat{v}_{2}(\xi,m^{\hat v}_{t}))}{\partial v_{2}}\partial_{m_{2}}\hat{v}_{2}(\xi,m^{\hat v}_{t})(x)
			\right) \right\} 
		m^{\hat v}_{2}(\xi,t)d\xi
	\right]dxdt
	\\
	&
	+\int_{0}^{T}\int_{\RR^{n}}\dfrac{\partial f_{1}(x,m^{\hat v}_{t},\hat{v}_{1}(x,m^{\hat v}_{t}))}{\partial v_{1}}\tilde{v}_{1}(x,m^{\hat v}_{t})m^{\hat v}_{1}(x,t)dxdt
	\\
	&
	+\int_{\RR^{n}}\tilde{m}_{1}(x,T)u^{\hat v}_{1}(x,T)dx+\int_{\RR^{n}}\tilde{m}_{-1}(x,T)u^{\hat v}_{-1}(x,T)dx \, .
\end{align*}
Integrating by parts and using the equations of $\tilde{m}_{1},\tilde{m}_{-1}$,
see (\ref{eq:NMFC-tildem1}), (\ref{eq:NMFC-tildem-1}), we obtain 

\begin{align}
	&\dfrac{d}{d\theta}\cJ^{NMFC}_{1}(\hat{v}_{1}+\theta\tilde{v}_{1})|_{\theta=0}
	\label{eq:2-9}
	\\
	=
	&\int_{0}^{T}\int_{\RR^{n}}\left[ \dfrac{\partial f_{1}(x,m^{\hat v}_{t},\hat{v}_{1}(x,m^{\hat v}_{t}))}{\partial v_{1}}
	+Du^{\hat v}_{1}(x,t).\dfrac{\partial g_{1}(x,m^{\hat v}_{t},\hat{v}_{1}(x,m^{\hat v}_{t}))}{\partial v_{1}} \right] \tilde{v}_{1}(x,m^{\hat v}_{t})m^{\hat v}_{1}(x,t)dxdt \, .
	\notag
\end{align}
Since this quantity must be $0$ for any possible $\tilde{v}_{1}(x,m)$,
we must have 
\begin{equation}
	\dfrac{\partial f_{1}(x,m^{\hat v}_{t},\hat{v}_{1}(x,m^{\hat v}_{t}))}{\partial v_{1}}
	+Du^{\hat v}_{1}(x,t).\dfrac{\partial g_{1}(x,m^{\hat v}_{t},\hat{v}_{1}(x,m^{\hat v}_{t}))}{\partial v_{1}}
	=0 \, .
	\label{eq:NMFC-1-optcond}
\end{equation}
So, using the Hamiltonian notation introduced in~\eqref{eq:defH-CMFC}, we have 
\begin{align}
	f_{1}(x,m^{\hat v}_{t},\hat{v}_{1}(x,m^{\hat v}_{t}))
	+Du^{\hat v}_{1}(x,t).g_{1}(x,m^{\hat v}_{t},\hat{v}_{1}(x,m^{\hat v}_{t}))
	= 
	H_{1}(x,m^{\hat v}_{t},Du^{\hat v}_{1}(x,t)) \, ,
	\label{eq:NMFC-f1g1H1}
\end{align}
and
\begin{align}
	&
	\int_{\RR^{n}} \left[ \partial_{m_{1}}f_{1}(\xi,m^{\hat v}_{t},\hat{v}_{1}(\xi,m^{\hat v}_{t}))(x)
		+\dfrac{\partial f_{1}(\xi,m^{\hat v}_{t},\hat{v}_{1}(\xi,m^{\hat v}_{t}))}{\partial v_{1}}\partial_{m_{1}}\hat{v}_{1}(\xi,m^{\hat v}_{t})(x) \right] 
		m^{\hat v}_{1}(\xi,t)d\xi
	\notag
	\\
	&+\int_{\RR^{n}} \left[ Du^{\hat v}_{1}(\xi,t).\left(
		\partial_{m_{1}}g_{1}(\xi,m^{\hat v}_{t},\hat{v}_{1}(\xi,m^{\hat v}_{t}))(x)
		+ \dfrac{\partial g_{1}(\xi,m^{\hat v}_{t},\hat{v}_{1}(\xi,m^{\hat v}_{t}))}{\partial v_{1}}\partial_{m_{1}}\hat{v}_{1}(\xi,m^{\hat v}_{t})(x)
		\right) \right] m^{\hat v}_{1}(\xi,t)d\xi
	\notag
	\\
	= \, 
	& \int_{\RR^{n}}\partial_{m_{1}}H_{1}(\xi,m^{\hat v}_{t},Du^{\hat v}_{1}(\xi,t))(x)m^{\hat v}_{1}(\xi,t)d\xi \, ,
	\label{eq:NMFC-dm1H1}
\end{align}
and
\begin{align}
	& \int_{\RR^{n}} \left[ \partial_{m_{2}}f_{1}(\xi,m^{\hat v}_{t},\hat{v}_{1}(\xi,m^{\hat v}_{t}))(x)+\dfrac{\partial f_{1}(\xi,m^{\hat v}_{t},\hat{v}_{1}(\xi,m^{\hat v}_{t}))}{\partial v_{1}}\partial_{m_{2}}\hat{v}_{1}(\xi,m^{\hat v}_{t})(x) \right] m^{\hat v}_{1}(\xi,t)d\xi
	\notag
	\\
	& +\int_{\RR^{n}} \left[ Du^{\hat v}_{1}(\xi,t).\left(
		\partial_{m_{2}}g_{1}(\xi,m^{\hat v}_{t},\hat{v}_{1}(\xi,m^{\hat v}_{t}))(x)
		+ \dfrac{\partial g_{1}(\xi,m^{\hat v}_{t},\hat{v}_{1}(\xi,m^{\hat v}_{t}))}{\partial v_{1}}\partial_{m_{2}}\hat{v}_{1}(\xi,m^{\hat v}_{t})(x)
	\right) \right] m^{\hat v}_{1}(\xi,t)d\xi
	\notag
	\\
	= \, 
	& \int_{\RR^{n}}\partial_{m_{2}}H_{1}(\xi,m^{\hat v}_{t},Du^{\hat v}_{1}(\xi,t))(x)m^{\hat v}_{1}(\xi,t)d\xi \, .
	\label{eq:NMFC-dm2H1}
\end{align}

We can write a similar necessary condition for player $2$ by introducing $u_2$ and $u_{-2}$. 
Unfortunately, this is not sufficient to write a self-contained system
for $u_{1},$ $u_{-1},$ $u_{2},$ $u_{-2},$ $m_{1},$ and $m_{2}$. Unlike the CMFC case, we
need to know explicitly the feedback $(x,m) \mapsto \hat{v}_{2}(x,m)$ to express
its functional derivatives $\partial_{m_{1}} \hat{v}_{2}(x,m)(\xi),$ $\partial_{m_{2}}\hat{v}_{2}(x,m)(\xi).$
This is a major difference with respect to the CMFC case, in which these
terms did not appear. For this reason, in the NMFC problem we cannot avoid the master equations in general.

\begin{remark}
\label{rem:NMFC-vanish}
In some special cases, it is not necessary to introduce the master equations; see section~\ref{sec:NMFC:specialcases} below.
\end{remark}

% subsec-subsec-subsec-subsec-subsec-subsec-subsec-subsec-subsec %
% subsec-subsec-subsec-subsec-subsec-subsec-subsec-subsec-subsec %
\subsection{Master equations}

We look for functions $(x,m,t) \mapsto U_{1}(x,m,t)$, $U_{-1}(x,m,t)$ such that
$u^{\hat v}_{1}(x,t)=U_{1}(x,m^{\hat v}_{t},t),$ $u^{\hat v}_{-1}(x,t)=U_{-1}(x,m^{\hat v}_{t},t)$ when $m^{\hat v}$ solves~\eqref{FPCeqGen} controlled by $\hat v$.
From (\ref{eq:NMFC-1-optcond}) we have, recalling the definition of $\hat{\bfv}_{1}$,
see (\ref{eq:defH-CMFC}),
\[
	\hat{v}_{1}(x,m^{\hat v}_{t})
	=
	\hat{\bfv}_{1}(x,m^{\hat v}_{t},Du^{\hat v}_{1}(x,t))
	=
	\hat{\bfv}_{1}(x,m^{\hat v}_{t},DU_{1}(x,m^{\hat v}_{t},t)) \, .
\] 
 We then rewrite~\eqref{eq:NMFC-f1g1H1}
as 
\begin{equation}
	f_{1}(x,m^{\hat v}_t,\hat{v}_{1}(x,m^{\hat v}_t))
	+Du^{\hat v}_{1}(x,t).g_{1}(x,m^{\hat v}_t,\hat{v}_{1}(x,m^{\hat v}_t))
	=
	H_{1}(x,m^{\hat v}_t,DU_{1}(x,m^{\hat v}_t,t)) \, ,
	\label{eq:2-13}
\end{equation}
and the right-hand sides of~\eqref{eq:NMFC-dm1H1} and~\eqref{eq:NMFC-dm2H1} become respectively
\begin{align}
	&\int_{\RR^{n}}\partial_{m_{1}}H_{1}(\xi,m^{\hat v}_t,Du^{\hat v}_{1}(\xi,t))(x)m^{\hat v}_{1}(\xi,t)d\xi 
	=
	\int_{\RR^{n}}\partial_{m_{1}}H_{1}(\xi,m^{\hat v}_t,DU_{1}(\xi,m^{\hat v}_t,t))(x)m^{\hat v}_{1}(\xi,t)d\xi \, ,
	\label{eq:2-14}
	\\
	& \int_{\RR^{n}}\partial_{m_{2}}H_{1}(\xi,m^{\hat v}_t,Du^{\hat v}_{1}(\xi,t))(x)m^{\hat v}_{1}(\xi,t)d\xi
	=
	\int_{\RR^{n}}\partial_{m_{2}}H_{1}(\xi,m^{\hat v}_t,DU_{1}(\xi,m^{\hat v}_t,t))(x)m^{\hat v}_{1}(\xi,t)d\xi \, .
	\label{eq:2-15}
\end{align}
In addition, using~\eqref{eq:NMFC-1-optcond}, we have the relation 
\begin{equation}
	g_{1}(x,m^{\hat v}_t,\hat{v}_{1}(x,m^{\hat v}_t))
	=
	\dfrac{\partial H_{1}}{\partial q_{1}}(x,m^{\hat v}_t,D_{x}U_{1}(x,m^{\hat v}_t,t)) \, .
	\label{eq:2-16}
\end{equation}
We can then take the functional derivatives of both sides, with
respect to $m_{1}$ and $m_{2}$ to obtain 
\begin{align}
	&\partial_{m_{1}}g_{1}(x,m^{\hat v}_t,\hat{v}_{1}(x,m^{\hat v}_t))(\xi)
	+\dfrac{\partial g_{1}}{\partial v_{1}}(x,m^{\hat v}_t,\hat{v}_{1}(x,m^{\hat v}_t)) \partial_{m_{1}}\hat{v}_{1}(x,m^{\hat v}_t)(\xi)
	\notag
	\\
	=\,
	&
	\partial_{m_{1}}\dfrac{\partial H_{1}}{\partial q_{1}}(x,m^{\hat v}_t,D_{x}U_{1}(x,m^{\hat v}_t,t))(\xi)
	+\dfrac{\partial^{2}H_{1}}{\partial q_{1}^{2}}(x,m^{\hat v}_t,D_{x}U_{1}(x,m^{\hat v}_t,t))\partial_{m_{1}}D_{x}U_{1}(x,m^{\hat v}_t,t)(\xi) \, ,
	\label{eq:NMFC-dm1g1}
\end{align}
and
\begin{align}
	& \partial_{m_{2}}g_{1}(x,m^{\hat v}_t,\hat{v}_{1}(x,m^{\hat v}_t))(\xi)
	+\dfrac{\partial g_{1}}{\partial v_{1}}(x,m^{\hat v}_t,\hat{v}_{1}(x,m^{\hat v}_t))\partial_{m_{2}}\hat{v}_{1}(x,m^{\hat v}_t)(\xi)
	\notag
	\\
	= \,
	& 
	\partial_{m_{2}}\dfrac{\partial H_{1}}{\partial q_{1}}(x,m^{\hat v}_t,D_{x}U_{1}(x,m^{\hat v}_t,t))(\xi)
	+ \dfrac{\partial^{2}H_{1}}{\partial q_{1}^{2}}(x,m^{\hat v}_t,D_{x}U_{1}(x,m^{\hat v}_t,t))\partial_{m_{2}}D_{x}U_{1}(x,m^{\hat v}_t,t)(\xi) \, .
	\label{eq:NMFC-dm2g1}
\end{align}

In equations~\eqref{eq:NMFC-u1-hat}, \eqref{eq:NMFC-u-1-hat} we do not need directly~\eqref{eq:NMFC-dm1g1}, \eqref{eq:NMFC-dm2g1} but their equivalent for the feedback
$(x,m) \mapsto \hat{v}_{2}(x,m)$, which involve functions $U_{2}$, $U_{-2}$ such that $u^{\hat v}_{2}(x,t) = U_{2}(x, m^{\hat v}_t, t)$, $u^{\hat v}_{-2}(x,t) = U_{-2}(x, m^{\hat v}_t, t)$. 
Finally, we can write~\eqref{eq:NMFC-u1-hat}, \eqref{eq:NMFC-u-1-hat} as 
\begin{align}
	&-\dfrac{\partial u_{1}}{\partial t}(x,t) + A_{1}u_{1}(x,t)
	\label{eq:NMFC-sys-HJB11}
	\\
	= \, 
	& H_{1}(x,m_{t},Du_{1}(x,t))+\int_{\RR^{n}}\partial_{m_{1}}H_{1}(\xi,m_{t},Du_{1}(\xi,t))(x)m_{1}(\xi,t)d\xi
	\notag
	\\
	& +\int_{\RR^{n}}Du_{-1}(\xi,t). 
	\left[ 
		\partial_{m_{1}}\dfrac{\partial H_{2}}{\partial q_{2}}(\xi,m_{t},D_{\xi}u_{2}(\xi,t))(x)
		\vphantom{\dfrac{\partial^{2}H_{2}}{\partial q_{2}^{2}}}
		\right.
	\notag
	\\
	&\qquad\qquad \left.
		+\dfrac{\partial^{2}H_{2}}{\partial q_{2}^{2}}(\xi,m_{t},D_{x}u_{2}(\xi,t))\partial_{m_{1}}D_{\xi}U_{2}(\xi,m_{t},t)(x)
	\right] m_{2}(\xi,t)d\xi \, ,
	\notag
\end{align}
 and 
\begin{align}
	&-\dfrac{\partial u_{-1}}{\partial t}(x,t) + A_{2}u_{-1} (x,t)
	\label{eq:NMFC-sys-HJB1-1}
	\\
	= \, 
	& D_{x}u_{-1}(x,t).\dfrac{\partial H_{2}}{\partial q_{2}}(x,m_{t},D_{x}u_{2}(x,t))+\int_{\RR^{n}}\partial_{m_{2}}H_{1}(\xi,m_{t},Du_{1}(\xi,t))(x)m_{1}(\xi,t)d\xi
	\notag
	\\
	& +\int_{\RR^{n}}Du_{-1}(\xi,t).
		\left[
			\partial_{m_{2}}\dfrac{\partial H_{2}}{\partial q_{2}}(\xi,m_{t},D_{\xi}u_{2}(\xi,t))(x)
		\vphantom{\dfrac{\partial^{2}H_{2}}{\partial q_{2}^{2}}}
		\right.
	\notag
	\\
	&\qquad\qquad \left.
			+\dfrac{\partial^{2}H_{2}}{\partial q_{2}^{2}}(\xi,m_{t},D_{\xi}u_{2}(\xi,t))\partial_{m_{2}}D_{\xi}U_{2}(\xi,m_{t},t)(x)
		\right] m_{2}(\xi,t)d\xi \, ,
	\notag
\end{align}
with terminal conditions 
\begin{align}
	u_{1}(x,T)
	&=
	h_{1}(x,m_{T})+\int_{\RR^{n}}\partial_{m_{1}}h_{1}(\xi,m_{T})(x)m_{1}(\xi,T)d\xi
	\notag
	%\label{eq:2-21}
	\\
	u_{-1}(x,T)
	&=
	\int_{\RR^{n}}\partial_{m_{2}}h_{1}(\xi,m_{T})(x)m_{1}(\xi,T)d\xi \, .
	\notag
\end{align}
Similarly, we associate equations for $u_{2}$ and $u_{-2}$ as follows

\begin{align}
	&-\dfrac{\partial u_{2}}{\partial t}(x,t) + A_{2}u_{2}(x,t)
	\label{eq:NMFC-sys-HJB22}
	\\
	= \,
	& H_{2}(x,m_{t},Du_{2}(x,t))+\int_{\RR^{n}}\partial_{m_{2}}H_{2}(\xi,m_{t},Du_{2}(\xi,t))(x)m_{2}(\xi,t)d\xi
	\notag
	\\
	& +\int_{\RR^{n}}Du_{-2}(\xi,t).
		\left[
			\partial_{m_{2}}\dfrac{\partial H_{1}}{\partial q_{1}}(\xi,m_{t},D_{\xi}u_{1}(\xi,t))(x)
		\vphantom{\dfrac{\partial^{2}H_{2}}{\partial q_{2}^{2}}}
		\right.
	\notag
	\\
	&\qquad\qquad \left.
		+\dfrac{\partial^{2}H_{1}}{\partial q_{1}^{2}}(\xi,m_{t},D_{x}u_{1}(\xi,t))\partial_{m_{2}}D_{\xi}U_{1}(\xi,m_{t},t)(x) 
	\right] m_{1}(\xi,t)d\xi \, ,
	\notag
\end{align}
and 
\begin{align}
	&-\dfrac{\partial u_{-2}}{\partial t}(x,t) + A_{1}u_{-2}(x,t)
	\label{eq:NMFC-sys-HJB2-2}
	\\
	= \, 
	& D_{x}u_{-2}(x,t).\dfrac{\partial H_{1}}{\partial q_{1}}(x,m_{t},D_{x}u_{1}(x,t))
	+\int_{\RR^{n}}\partial_{m_{1}}H_{2}(\xi,m_{t},Du_{2}(\xi,t))(x)m_{2}(\xi,t)d\xi
	\notag
	\\
	& +\int_{\RR^{n}}Du_{-2}(\xi,t).
	\left[ 
		\partial_{m_{1}}\dfrac{\partial H_{1}}{\partial q_{1}}(\xi,m_{t},D_{\xi}u_{1}(\xi,t))(x)
		\vphantom{\dfrac{\partial^{2}H_{2}}{\partial q_{2}^{2}}}
		\right.
	\notag
	\\
	&\qquad\qquad \left.
		+\dfrac{\partial^{2}H_{1}}{\partial q_{1}^{2}}(\xi,m_{t},D_{\xi}u_{1}(\xi,t))\partial_{m_{1}}D_{\xi}U_{1}(\xi,m_{t},t)(x)
	\right] m_{1}(\xi,t)d\xi \, ,
	\notag
\end{align}
with the terminal conditions
\begin{align}
	u_{2}(x,T)
	&=
	h_{2}(x,m_{T})+\int_{\RR^{n}}\partial_{m_{2}}h_{2}(\xi,m_{T})(x)m_{2}(\xi,T)d\xi \, ,
	\notag
	%\label{eq:2-24}
	\\
	u_{-2}(x,T)
	&=
	\int_{\RR^{n}}\partial_{m_{1}}h_{2}(\xi,m_{T})(x)m_{2}(\xi,T)d\xi \, ,
	\notag
\end{align}
 and the FP equations are 
\begin{align}
	& \dfrac{\partial m_{1}}{\partial t}+A_{1}^{*}m_{1}+\diver_{x}\left( \dfrac{\partial H_{1}}{\partial q_{1}}(x, m_{t},D_{x}u_{1}(x,t))m_{1} \right)
	=
	0
	%\notag
	\label{eq:NMFC-sys-FP1}
	\\
	& \dfrac{\partial m_{2}}{\partial t}+A_{2}^{*}m_{2}+\diver_{x}\left(\dfrac{\partial H_{2}}{\partial q_{2}}(x, m_{t},D_{x}u_{2}(x,t))m_{2} \right)
	=
	0
	%\notag
	\label{eq:NMFC-sys-FP2}
	\\
	& m_{1}(x,0)=\rho_{10}(x) \, ,\qquad m_{2}(x,0) = \rho_{20}(x) \, .
	\notag
\end{align}

A crucial remark is that the system of $6$ equations composed of~\eqref{eq:NMFC-sys-HJB11}, \eqref{eq:NMFC-sys-HJB1-1}, \eqref{eq:NMFC-sys-HJB22}, \eqref{eq:NMFC-sys-HJB2-2}, \eqref{eq:NMFC-sys-FP1}, and~\eqref{eq:NMFC-sys-FP2}, is \textbf{not self-contained} in general since we cannot express
$\partial_{m_{1}}D_{\xi}U_{1}(\xi,m_{t},t)(x),$ $\partial_{m_{2}}D_{\xi}U_{1}(\xi,m_{t},t)(x),$ $\partial_{m_{2}}D_{\xi}U_{2}(\xi,m_{t},t)(x),$ and $\partial_{m_{1}}D_{\xi}U_{2}(\xi,m_{t},t)(x)$
in terms of $u_{1},$ $u_{-1},$ $u_{2},$ and $u_{-2}.$ 
So we need to write the
system of equations for $U_{1},$ $U_{-1},$ $U_{2},$ and $U_{-2}$.
We obtain 
\begin{align}
	&-\dfrac{\partial U_{1}}{\partial t}(x,m,t) 
	+ A_{1}U_{1}(x,m,t) 
	\label{eq:2-26}
	\\
	= \, 
	&-\sum_{j=1}^2\int_{\RR^{n}}A_{j\xi}\partial_{m_{j}}U_{1}(\xi,m,t)(x)m_{j}(\xi)d\xi 
	\notag
	\\
	& +\sum_{j=1}^2 \int_{\RR^{n}}\partial_{m_{j}}D_{\xi}U_{1}(\xi,m,t)(x).\dfrac{\partial H_{j}}{\partial q_{j}}(\xi,m,D_{\xi}U_{j}(\xi,m,t))m_{j}(\xi)d\xi
	\notag
	\\
	& +H_{1}(x,m,D_{x}U_{1}(x,m,t))+\int_{\RR^{n}}\partial_{m_{1}}H_{1}(\xi,m,DU_{1}(\xi,m,t))(x)m_{1}(\xi)d\xi
	\notag
	\\
	& +\int_{\RR^{n}}DU_{-1}(\xi,m,t).
	\left[
		\partial_{m_{1}}\dfrac{\partial H_{2}}{\partial q_{2}}(\xi,m,D_{\xi}U_{2}(\xi,m,t))(x)
	\right.
	\notag
	\\
	&\qquad\qquad
	\left. 
		+\dfrac{\partial^{2}H_{2}}{\partial q_{2}^{2}}(\xi,m,D_{\xi}U_{2}(\xi,m,t))\partial_{m_{1}}D_{\xi}U_{2}(\xi,m,t)(x)
	\right] m_{2}(\xi)d\xi \, ,
	\notag
\end{align}
and
\begin{align}
	&-\dfrac{\partial U_{-1}}{\partial t}(x,m,t) + A_{2}U_{-1}(x,m,t) 
	\label{eq:2-27}
	\\
	= \, 
	& -\sum_{j=1}^2\int_{\RR^{n}}A_{j\xi}\partial_{m_{j}}U_{-1}(\xi,m,t)(x)m_{j}(\xi)d\xi  
	\notag
	\\
	& +\sum_{j=1}^2 \int_{\RR^{n}}\partial_{m_{j}}D_{\xi}U_{-1}(\xi,m,t)(x).\dfrac{\partial H_{j}}{\partial q_{j}}(\xi,m,D_{\xi}U_{j}(\xi,m,t))m_{j}(\xi)d\xi
	\notag
	\\
	& +D_{x}U_{-1}(x,m,t).\dfrac{\partial H_{2}}{\partial q_{2}}(x,m,D_{x}U_{2}(x,m,t))+\int_{\RR^{n}}\partial_{m_{2}}H_{1}(\xi,m,DU_{1}(\xi,m,t))(x)m_{1}(\xi)d\xi
	\notag
	\\
	& +\int_{\RR^{n}}DU_{-1}(\xi,m,t).
	\left[
		\partial_{m_{2}}\dfrac{\partial H_{2}}{\partial q_{2}}(\xi,m,D_{\xi}U_{2}(\xi,m,t))(x)
	\right.
	\notag
	\\
	& 
	\qquad\qquad 
	\left.
		+\dfrac{\partial^{2}H_{2}}{\partial q_{2}^{2}}(\xi,m,D_{\xi}U_{2}(\xi,m,t))\partial_{m_{2}}D_{\xi}U_{2}(\xi,m,t)(x)
	\right] m_{2}(\xi)d\xi \, ,
	\notag
\end{align}
with terminal conditions
\begin{align}
	U_{1}(x,m,T)
	& =
	h_{1}(x,m)+\int_{\RR^{n}}\partial_{m_{1}}h_{1}(\xi,m)(x)m_{1}(\xi)d\xi
	\notag
	%\label{eq:270}
	\\
	U_{-1}(x,m,T)
	&=
	\int_{\RR^{n}}\partial_{m_{2}}h_{1}(\xi,m)(x)m_{1}(\xi)d\xi \, .
	\notag
\end{align}

Moreover, we have two additional equations for $U_{2},U_{-2}$ 
\begin{align}
	&-\dfrac{\partial U_{2}}{\partial t}(x,m,t)  + A_{2}U_{2}(x,m,t) 
	\label{eq:2-28}
	\\
	= \,
	&
	 -\sum_{j=1}^2\int_{\RR^{n}}A_{j\xi}\partial_{m_{j}}U_{2}(\xi,m,t)(x)m_{j}(\xi)d\xi  
	\notag
	\\
	& +\sum_{j=1}^2 \int_{\RR^{n}}\partial_{m_{j}}D_{\xi}U_{2}(\xi,m,t)(x).\dfrac{\partial H_{j}}{\partial q_{j}}(\xi,m,D_{\xi}U_{j}(\xi,m,t))m_{j}(\xi)d\xi
	\notag
	\\
	& +H_{2}(x,m,D_{x}U_{2}(x,m,t))+\int_{\RR^{n}}\partial_{m_{2}}H_{2}(\xi,m,D_{\xi}U_{2}(\xi,m,t))(x)m_{2}(\xi,t)d\xi
	\notag
	\\
	& +\int_{\RR^{n}}DU_{-2}(\xi,m,t).
	\left[
		\partial_{m_{2}}\dfrac{\partial H_{1}}{\partial q_{1}}(\xi,m,D_{\xi}U_{1}(\xi,m,t))(x)
	\right.
	\notag
	\\
	& \qquad\qquad
	\left.
		+\dfrac{\partial^{2}H_{1}}{\partial q_{1}^{2}}(\xi,m,D_{\xi}U_{1}(\xi,m,t))\partial_{m_{2}}D_{\xi}U_{1}(\xi,m,t)(x)
	\right] m_{1}(\xi)d\xi \, ,
	\notag
\end{align}
 and
\begin{align}
	&-\dfrac{\partial U_{-2}}{\partial t}(x,m,t) + A_{1}U_{-2}(x,m,t) 
	\label{eq:2-29}
	\\
	= \,
	&  -\sum_{j=1}^2\int_{\RR^{n}}A_{j\xi}\partial_{m_{j}}U_{-2}(\xi,m,t)(x)m_{j}(\xi)d\xi  
	\notag
	\\
	& +\sum_{j=1}^2 \int_{\RR^{n}}\partial_{m_{j}}D_{\xi}U_{-2}(\xi,m,t)(x).\dfrac{\partial H_{j}}{\partial q_{j}}(\xi,m,D_{\xi}U_{j}(\xi,m,t))m_{j}(\xi)d\xi
	\notag
	\\
	& +D_{x}U_{-2}(x,m,t).\dfrac{\partial H_{1}}{\partial q_{1}}(x,m,D_{x}U_{1}(x,m,t))+\int_{\RR^{n}}\partial_{m_{1}}H_{2}(\xi,m,DU_{2}(\xi,m,t))(x)m_{2}(\xi)d\xi
	\notag
	\\
	& +\int_{\RR^{n}}DU_{-2}(\xi,m,t).
	\left[
		\partial_{m_{1}}\dfrac{\partial H_{1}}{\partial q_{1}}(\xi,m,D_{\xi}U_{1}(\xi,m,t))(x)
	\right.
	\notag
	\\
	&\qquad\qquad
	\left.
		+\dfrac{\partial^{2}H_{1}}{\partial q_{1}^{2}}(\xi,m,D_{\xi}U_{1}(\xi,m,t))\partial_{m_{1}}D_{\xi}U_{1}(\xi,m,t)(x)
	\right] m_{1}(\xi)d\xi \, ,
	\notag
\end{align}
with terminal conditions 
\begin{align}
	U_{2}(x,m,T)
	&=
	h_{2}(x,m)+\int_{\RR^{n}}\partial_{m_{2}}h_{2}(\xi,m)(x)m_{2}(\xi)d\xi \, ,
	\notag
	%\label{eq:2-290}
	\\
	U_{-2}(x,m,T)
	&=
	\int_{\RR^{n}}\partial_{m_{1}}h_{2}(\xi,m)(x)m_{2}(\xi)d\xi \, .
	\notag
\end{align}
The system of four master equations~\eqref{eq:2-26}, \eqref{eq:2-27}, \eqref{eq:2-28},
\eqref{eq:2-29}  for the functions $U_{1},$ $U_{-1},$ $U_{2},$ $U_{-2}$ is self-contained. 

Recall that in our notation, the subscript ``$-1$'' does not mean ``$2$''. Here $U_{-1}$ stems from the problem of player $1$ and reflects the impact of the variation of the state of the opponent. This intuition is made more precise in the next subsection.

% subsec-subsec-subsec-subsec-subsec-subsec-subsec-subsec-subsec %
% subsec-subsec-subsec-subsec-subsec-subsec-subsec-subsec-subsec %

\subsection{Bellman system}

We can check that 
\begin{align}
	& U_{1}(x,m,t)
	=
	\partial_{m_{1}}V_{1}(m,t)(x) \, ,
	\qquad 
	U_{-1}(x,m,t)
	=
	\partial_{m_{2}}V_{1}(m,t)(x) \, ,
	\label{eq:NMFC-U-dmV1}
	%\label{eq:2-30}
	\\
	& U_{2}(x,m,t)
	=
	\partial_{m_{2}}V_{2}(m,t)(x) \, ,
	\qquad
	U_{-2}(x,m,t)
	=
	\partial_{m_{1}}V_{2}(m,t)(x) \, ,
	\label{eq:NMFC-U-dmV2}
\end{align}
where $(m,t) \mapsto V_{1}(m,t), V_{2}(m,t)$ solve of the Bellman system

\begin{align}
	& 0
	= \,
	\dfrac{\partial V_{1}}{\partial t}(m,t)
		-\int_{\RR^{n}}A_{1x}\partial_{m_{1}}V_{1}(m,t)(x)m_{1}(x)dx
		-\int_{\RR^{n}}A_{2x}\partial_{m_{2}}V_{1}(m,t)(x)m_{2}(x)dx
	\label{eq:2-31}
	\\
	& \qquad+\int_{\RR^{n}}H_{1}(x,m,D_{x}\partial_{m_{1}}V_{1}(m,t)(x))m_{1}(x)dx
	\notag
	\\
	& \qquad +\int_{\RR^{n}}D_{x}\partial_{m_{2}}V_{1}(m,t)(x).\dfrac{\partial H_{2}}{\partial q_{2}}(x,m,D_{x}\partial_{m_{2}}V_{2}(m,t)(x))m_{2}(x)dx \, ,
	\notag
%	\\
%	& V_{1}(m,T)
%		= \int_{\RR^{n}}h_{1}(x,m)m_{1}(x)dx \, ,
\end{align}
with terminal condition
$$
	V_{1}(m,T)
		= \int_{\RR^{n}}h_{1}(x,m)m_{1}(x)dx \, ,
$$
and
\begin{align}
	& 0 = \,
	\dfrac{\partial V_{2}}{\partial t}(m,t)
		-\int_{\RR^{n}}A_{1x}\partial_{m_{1}}V_{2}(m,t)(x)m_{1}(x)dx
		-\int_{\RR^{n}}A_{2x}\partial_{m_{2}}V_{2}(m,t)(x)m_{2}(x)dx
	\label{eq:2-32}
	\\
	& \qquad +\int_{\RR^{n}}H_{2}(x,m,D_{x}\partial_{m_{2}}V_{2}(m,t)(x))m_{2}(x)dx
	\notag
	\\
	& \qquad +\int_{\RR^{n}}D_{x}\partial_{m_{1}}V_{2}(m,t)(x).\dfrac{\partial H_{1}}{\partial q_{1}}(x,m,D_{x}\partial_{m_{1}}V_{1}(m,t)(x))m_{1}(x)dx \, ,
	\notag
%	\\
%	& V_{2}(m,T)
%		= \int_{\RR^{n}}h_{2}(x,m)m_{2}(x)dx \, .
\end{align}
 with terminal condition
 $$
 	 V_{2}(m,T)
		= \int_{\RR^{n}}h_{2}(x,m)m_{2}(x)dx \, .
$$
\begin{remark}
To recover the equations (\ref{eq:2-26}) to (\ref{eq:2-29}) by taking
the functional derivatives of Bellman equations, one needs to use
the following symmetry properties stemming from~\eqref{eq:NMFC-U-dmV1}--\eqref{eq:NMFC-U-dmV2} 
\[
	\partial_{m_{2}}U_{1}(x,m,t)(\xi)=\partial_{m_{1}}U_{-1}(\xi,m,t)(x) \, ,
	\qquad
	\partial_{m_{1}}U_{2}(x,m,t)(\xi)=\partial_{m_{2}}U_{-1}(\xi,m,t)(x) \, .
\]
\end{remark}

\begin{remark}
	After completion of this work, it has been brought to our attention that similar problems have been studied recently under the term \textbf{mean field type games}, see e.g.~\cite{Math-02-00706} for the general setting and~\cite{ElectronEng-01-00018} for a collection of applications to engineering. However, to the best of our knowledge, our work is the first, at this level of generality, to provide a comprehensive framework and to focus on the necessary conditions of optimality formulated in terms of PDEs.
\end{remark}

% secsecsecsecsecsecsecsecsecsecsecsecsecsecsecsecsecsecsecsecsecsecsecsec %
% secsecsecsecsecsecsecsecsecsecsecsecsecsecsecsecsecsecsecsecsecsecsecsec %
% secsecsecsecsecsecsecsecsecsecsecsecsecsecsecsecsecsecsecsecsecsecsecsec %

\section{Mean Field Game Problems}
\label{sec:MFG}

In this section, we consider the interaction of two populations when the agents of each population are rational (i.e., try to minimize an individual cost and anticipate the rationality of other the players) and we look for non-cooperative equilibrium using a mean field game approach. For a single population, the mean field game viewpoint focuses on Nash equilibria among the population and, as such, the problem is defined through a fixed point procedure: first, given the distribution of the other players' states, an infinitesimal player finds her best response (that is, her optimal control) in order to minimize her cost; second, the distribution driven by the optimal control found in the first step should correspond to the distribution of the population. The reader is referred to e.g.~\cite{MR3134900} (Chapter 2) and~\cite{MR3045029} for more details. This idea can be extended to the case of several populations.
The PDE system is expressed e.g. in~\cite{MR3134900} (Chapter 8, page 68) and in~\cite{MR3752669} (Chapter 7, page 625) in the setting that has been the most usual in the literature so far. 
Analogously to the mean field control approach, we will distinguish between two types of problems.

% subsec-subsec-subsec-subsec-subsec-subsec-subsec-subsec-subsec %
% subsec-subsec-subsec-subsec-subsec-subsec-subsec-subsec-subsec %

\subsection{Nash Mean Field Game}
%\modif{
Let us start with a setting where each player compete with all the other players. We call it \textbf{Nash mean field game} (NMFG for short). 
In this problem each infinitesimal agent considers that the distributions of both populations are fixed and tries to minimize her own individual cost. Fixing the state of all the other players is translated, in the mean field limit, by the fact that the probabilities $m = (m_1,m_2)^*$ entering in the functions
$f_{i}(x,m,v_{i}),$ $g_{i}(x,m,v_{i}),$ $h_{i}(x,m)$ are considered
as fixed parameters. For this reason, we look for feedbacks depending on $x$ only instead of $(x,m)$, that is, the controls are functions $x \mapsto v_i(x)$.

\begin{problem}[NMFG]\label{pb:NMFG} Find $(\hat m, \hat v)$ satisfying the two conditions
\begin{enumerate}
	\item $\hat v = (\hat v_1, \hat v_2)^*$ is a Nash equilibrium for the functionals
	\begin{align*}
		&J^{NMFG}_{i,\hat m}(v_{1}, v_{2})
%		\notag
%		\\
		= \,
%		&\int_{0}^{T}\int_{\RR^{n}}f_{i}(x,\hat m_{t},v_{i}(x,\hat m_{t})) m^{v_{i},\hat v_{-i}, \hat m}_{i}(x,t) dxdt+\int_{\RR^{n}}h_{i}(x,\hat m_{T})m^{v_{i},\hat v_{-i}, \hat m}_{i}(x,T) dx
		&\int_{0}^{T}\int_{\RR^{n}}f_{i}(x,\hat m_{t},v_{i}(x)) m^{v, \hat m}_{i}(x,t) dxdt+\int_{\RR^{n}}h_{i}(x,\hat m_{T})m^{v, \hat m}_{i}(x,T) dx
		%\label{eq:2-1}
	\end{align*}
%where $m^{v_{i},\hat v_{-i}, \hat m}= (m^{v_{i},\hat v_{-i}, \hat m}_1,m^{v_{i},\hat v_{-i}, \hat m}_2)$ satisfies the PDEs~\eqref{FPCeqGen-MFG} controlled by $(v_{i},\hat v_{-i})$
where $m^{v, \hat m}= (m^{v, \hat m}_1,m^{v, \hat m}_2)^*$ satisfies
	\begin{equation*}%\label{FPCeqGen-MFG}
		\frac{\partial m_i}{\partial t}(x,t)+A_i^*m_i(x,t)+\diver_x\, \big(g_i(x,\hat m_t,v_i(x))m_i(x,t)\big) 
		=0 \, ,
	\end{equation*}
	with initial conditions $m_i(\cdot,0)=\rho_{i,0}$.
	\item $\hat m = m^{\hat v} = (m^{\hat v}_1,m^{\hat v}_2)^*$ is a  solution to~\eqref{FPCeqGen-MFG} controlled by $\hat v$. 
\end{enumerate}
\end{problem}
The first condition means that, for a given $\hat m$, for any $v$,
$$
	J^{NMFG}_{1,\hat m}(\hat v_{1}, \hat v_{2}) \leq J^{NMFG}_{1,\hat m}(v_{1}, \hat v_{2})
	\quad \hbox{ and } \quad
	J^{NMFG}_{2,\hat m}(\hat v_{1}, \hat v_{2}) \leq J^{NMFG}_{2,\hat m}(\hat v_{1}, v_{2}) \, .
$$

The problem of player $1$ is the following.
\begin{problem}[NMFG: Problem of player 1] 
\label{pb:NMFG-player1}
Find $\hat v_1$ minimizing
\begin{align*}
		&{\mathcal J}^{NMFG}_{1}(v_{1}) = J^{NMFG}_{1,\hat m}(v_{1}, \hat v_{2})
		\notag
		\\
		= \,
		&\int_{0}^{T}\int_{\RR^{n}} f_{1}\left(x, \hat m_{t},v_{1}\left(x\right)\right) m^{v_{1},\hat v_{2}, \hat m}_{1}(x,t) dxdt
%		\\
%		&\qquad 
		+\int_{\RR^{n}} h_{1}\left(x, \hat m_{T} \right)m^{v_{1},\hat v_{2}, \hat m}_{1}(x,T) dx
	\end{align*}
where $\left( m^{v_{1},\hat v_{2}, \hat m}_1, m^{v_{1},\hat v_{2}, \hat m}_2\right)$ solves
	\begin{align}
		&\frac{\partial m_1}{\partial t}(x,t)
		+A_1^*m_1(x,t)
		+\diver_x\, \Big(g_1\left(x,\hat m_{t},v_{1}\left(x\right)\right)m_1(x,t)\Big) 
		=0 \, ,
		\label{FPCeqGen-NMFG-THML-1-1}
		\\
		&\frac{\partial m_2}{\partial t}(x,t)
		+A_2^*m_2(x,t)
		+\diver_x\, \Big(g_2\left(x, \hat m_{t}, \hat v_{2}\left(x\right)\right)m_2(x,t)\Big) 
		=0 \, ,
		\label{FPCeqGen-NMFG-THML-1-2}
	\end{align}
	with initial conditions $m_i(\cdot,0)=\rho_{i,0}$.
\end{problem}
By comparing Problem~\ref{pb:NMFG-player1} and Problem~\ref{pb:NMFC-player1}, we see that in order to derive the necessary optimality conditions, we can reuse the computations done in the NMFC setting but now with a fixed parameter $\hat m = (\hat m_1, \hat m_2)^*$, and impose a posteriori that $\hat m = (\hat m_1, \hat m_2)^*$ should solve the FP equations. We shall only provide the equations and skip the proof to alleviate the presentation. This leads to the PDE system
\begin{align}
	&-\dfrac{\partial u_{i}}{\partial t}(x,t) + A_{i}u_{i}(x,t)
	=
	H_{i}(x,m_{t},Du_{i}(x,t))
	\label{eq:CMFG-HJBFPsys-HJB}
	%\label{eq:3-1}
	\\
	&\dfrac{\partial m_{i}}{\partial t}(x,t) + A_{i}^{*}m_{i}(x,t) + \diver_{x}\left(\dfrac{\partial H_{i}}{\partial q_{i}}(x,m_{t},Du_{i}(x,t))m_{i}(x,t)\right)
	=0
	\label{eq:CMFG-HJBFPsys-FP}
\end{align}
with terminal and initial conditions
$$
	u_{i}(x,T)=h_{i}(x,m_{T}) \, ,\qquad m_{i}(x,0)=\rho_{i,0}(x) \, .
$$
	This corresponds to the PDE system obtained in~\cite{MR3134900} (Chapter 8, page 68), for which the mean field problem had not been written explicitly but which has been shown to provide an approximate Nash equilibrium for a finite player game where all the players compete. We refer the interested reader to~\cite{MR3134900} and~\cite{MR3752669} for more details.

% subsec-subsec-subsec-subsec-subsec-subsec-subsec-subsec-subsec %
% subsec-subsec-subsec-subsec-subsec-subsec-subsec-subsec-subsec %

\subsection{Common Mean Field Game}

We now consider a different viewpoint where, in analogy with the CMFC setting, there is only one (common) cost functional. 
We introduce the following problem, that we call \textbf{common mean field game} (CMFG for short).
 Here again, the states of the populations are fixed so we look for feedbacks as functions of $x$ only instead of $(x,m)$.
\begin{problem}[CMFG]\label{pb:CMFG}
 Find $(\hat m, \hat v)$ satisfying the two conditions
 \begin{enumerate}
 \item $\hat v=(\hat v_1,\hat v_2)^*$ minimizes
	\begin{align}
		J^{CMFG}_{\hat m}(v_1,v_2)
		&=
		\sum_{i=1}^2  \left[ \int_0^T \int_{\RR^n} f_i(x,\hat m_t,v_i(x)) m^{v, \hat m}_i(x,t) dx \,dt + \int_{\RR^n} h_i(x,\hat m_T) m^{v, \hat m}_i(x,T) dx \right],
		\label{defJ-CMFG}
	\end{align}
	where $m^{v, \hat m}= (m^{v, \hat m}_1,m^{v, \hat m}_2)^*$ satisfies the following PDE system
	\begin{equation}\label{FPCeqGen-MFG}
		\frac{\partial m_i}{\partial t}(x,t)+A_i^*m_i(x,t)+\diver_x\, \big(g_i(x,\hat m_t,v_i(x))m_i(x,t)\big) 
		=0 \, ,
	\end{equation}
	with initial conditions $m_i(\cdot,0)=\rho_{i,0}$.
   \item
   $\hat m = m^{\hat v,\hat m}$ is a solution to~\eqref{FPCeqGen-MFG} controlled by $\hat v$. 
\end{enumerate}
\end{problem}

The control problem appearing in the first point above is solved ignoring the
parameters $\hat m_1, \hat m_2$. Eventually the value of these parameters is defined a
posteriori, by a fixed point argument, equaling these parameters to
the solution of the FP equations~\eqref{FPCeqGen-MFG}. We can use this viewpoint to apply
the CMFC considered in section~\ref{sec:CMFC}. Referring to the system~\eqref{eq:CMFC-HJBFPsys-HJB} of HJB-FP equations, 
 it turns out that we recover the equations~\eqref{eq:CMFG-HJBFPsys-HJB}--\eqref{eq:CMFG-HJBFPsys-FP} with the same terminal and initial conditions. Hence the necessary optimality conditions of the two mean field games (Problems~\ref{pb:NMFG} and~\ref{pb:CMFG}) have the same PDE system. In other words, we have two different interpretations for this PDE system. This similarity between the NMFG and the CMFG problems is not entirely surprising since, by looking at the definition of Problem~\ref{pb:CMFG}, one realizes that the problem can be split into two sub-problems (one minimization problem for each component of $v$) which are independent because $\hat m = (\hat m_1,\hat m_2)^*$ is fixed.

Although, as described above, one can reuse the computations done in the CMFC setting to derive the PDE system of  CMFG, the difference between the two PDE systems should be stressed: the HJB equation~\eqref{eq:CMFG-HJBFPsys-HJB} does not involve derivatives of the Hamiltonians with respect to $m_i$, $i=1,2,$ whereas in~\eqref{eq:CMFC-HJBFPsys-HJB} such derivatives do appear.

\begin{remark}
Studying rigorously the corresponding $N$-player game is an important question which is beyond the scope of the present work. However, at least at a heuristic level, the CMFG problem described above can be viewed as the mean field limit of a game with a finite number of players in at least two ways. One could imagine that the state of each player has two components, each with its own dynamics, and the player chooses a control for each component in order to minimize her global cost (which depends on both components). One could also consider a game with two populations of equal size, say $(X^i_1)_{i=1,\dots,N}$ and $(X^i_2)_{i=1,\dots,N}$, where the players work in pairs: $X^i_1$ and $X^i_2$ help each other and compete with $((X^j_k)_{k=1,2})_{j \neq i}$. This would be a game representing competition between couples composed of one player from each population: the players collaborate among each pair but compete at a global level.
\end{remark}

\begin{remark} 
The adjoint equations~\eqref{eq:CMFG-HJBFPsys-HJB} may also be deduced with an approach based on the master equation point of view. We omit the details.  
\end{remark}

% secsecsecsecsecsecsecsecsecsecsecsecsecsecsecsecsecsecsecsecsecsecsecsec %
% secsecsecsecsecsecsecsecsecsecsecsecsecsecsecsecsecsecsecsecsecsecsecsec %
% secsecsecsecsecsecsecsecsecsecsecsecsecsecsecsecsecsecsecsecsecsecsecsec %
\section{Examples}
\label{sec:examples}

%\modif{
% subsec-subsec-subsec-subsec-subsec-subsec-subsec-subsec-subsec %
% subsec-subsec-subsec-subsec-subsec-subsec-subsec-subsec-subsec %
\subsection{Special cases}
\label{sec:NMFC:specialcases}
Let us start with some situations in which the system of master equations is not needed and one can work instead with a system of PDEs in finite dimension. For future reference (see also the examples below), we report here two cases of interest for many applications. To the best of our knowledge, the examples studied in the literature so far fall in one of these cases.

\paragraph{Special case 1: } 
When deriving the necessary optimality conditions in the previous sections, we have looked for controls under the form of functions of both $x$ and $m = (m_1,m_2)^*$. This corresponds to a situation where each player observes her individual state together with both distributions. If, instead, one considers a more restricted information structure according to which the agents do not have access to the distributions, then one should look for controls under the form of functions of $x$ only, i.e. $x \mapsto v_i(x)$. In this case, the optimal controls, in particular, are not allowed to depend on $m$. Hence in~\eqref{eq:NMFC-u1-hat}--\eqref{eq:NMFC-u-1-hat} and in the analogous equations for $u_2,u_{-2}$, the terms $(\partial_{m_{i}}\hat{v}_{j})_{i=1,2,j=1,2}$ vanish. We can thus  write a self-contained PDE system for $u_{1},$ $u_{-1},$ $u_{2},$ $u_{-2},$ $m_{1},$ and $m_{2}$ in the NMFC setting. Indeed, the equations~\eqref{eq:NMFC-sys-HJB11}--\eqref{eq:NMFC-sys-HJB1-1} for $u_1, u_{-1}$ of NMFC become
\begin{align}
	&-\dfrac{\partial u_{1}}{\partial t}(x,t) + A_{1}u_{1}(x,t)
	\label{eq:NMFC-simplecase1-sys-HJB11}
	\\
	= \, 
	& H_{1}(x,m_{t},Du_{1}(x,t))+\int_{\RR^{n}}\partial_{m_{1}}H_{1}(\xi,m_{t},Du_{1}(\xi,t))(x)m_{1}(\xi,t)d\xi
	\notag
	\\
	& +\int_{\RR^{n}}Du_{-1}(\xi,t). 
%	\left[ 
		\partial_{m_{1}}\dfrac{\partial H_{2}}{\partial q_{2}}(\xi,m_{t},D_{\xi}u_{2}(\xi,t))(x)
		\vphantom{\dfrac{\partial^{2}H_{2}}{\partial q_{2}^{2}}}
%		\right.
%	\notag
%	\\
%	&\qquad\qquad \left.
%		+\dfrac{\partial^{2}H_{2}}{\partial q_{2}^{2}}(\xi,m_{t},D_{x}u_{2}(\xi,t))\partial_{m_{1}}D_{\xi}U_{2}(\xi,m_{t},t)(x)
%	\right] 
	m_{2}(\xi,t)d\xi \, ,
	\notag
\end{align}
 and 
\begin{align}
	&-\dfrac{\partial u_{-1}}{\partial t}(x,t) + A_{2}u_{-1} (x,t)
	\label{eq:NMFC-simplecase1-sys-HJB1-1}
	\\
	= \, 
	& D_{x}u_{-1}(x,t).\dfrac{\partial H_{2}}{\partial q_{2}}(x,m_{t},D_{x}u_{2}(x,t))+\int_{\RR^{n}}\partial_{m_{2}}H_{1}(\xi,m_{t},Du_{1}(\xi,t))(x)m_{1}(\xi,t)d\xi
	\notag
	\\
	& +\int_{\RR^{n}}Du_{-1}(\xi,t).
%		\left[
			\partial_{m_{2}}\dfrac{\partial H_{2}}{\partial q_{2}}(\xi,m_{t},D_{\xi}u_{2}(\xi,t))(x)
%		\vphantom{\dfrac{\partial^{2}H_{2}}{\partial q_{2}^{2}}}
%		\right.
%	\notag
%	\\
%	&\qquad\qquad \left.
%			+\dfrac{\partial^{2}H_{2}}{\partial q_{2}^{2}}(\xi,m_{t},D_{\xi}u_{2}(\xi,t))\partial_{m_{2}}D_{\xi}U_{2}(\xi,m_{t},t)(x)
%		\right] 
		m_{2}(\xi,t)d\xi \, ,
	\notag
\end{align}
with terminal conditions 
\begin{align}
	u_{1}(x,T)
	&=
	h_{1}(x,m_{T})+\int_{\RR^{n}}\partial_{m_{1}}h_{1}(\xi,m_{T})(x)m_{1}(\xi,T)d\xi
	\notag
	%\label{eq:2-21}
	\\
	u_{-1}(x,T)
	&=
	\int_{\RR^{n}}\partial_{m_{2}}h_{1}(\xi,m_{T})(x)m_{1}(\xi,T)d\xi \, .
	\notag
\end{align}
For $u_{2}$ and $u_{-2}$ similar equations hold, and the FP equations remain~\eqref{eq:NMFC-sys-FP1}--\eqref{eq:NMFC-sys-FP2}.

\paragraph{Special case 2: } If the controls are functions of $x$ only (as in the first special case above), and in addition $g_1$ does not depend upon $m_2$ and $g_2$ does not depend on $m_1$, then the unknowns $u_{-1}$ and $u_{-2}$ become superfluous and we can write a self-contained PDE system for $u_1,u_2,m_1,m_2$. Indeed, for NMFC, the equation~\eqref{eq:NMFC-simplecase1-sys-HJB11} for $u_1$ simplifies further and we obtain
\begin{align}
	&-\dfrac{\partial u_{1}}{\partial t}(x,t) + A_{1}u_{1}(x,t)
	\label{eq:NMFC-simplecase2-sys-HJB11}
	\\
	= \, 
	& H_{1}(x,m_{t},Du_{1}(x,t))+\int_{\RR^{n}}\partial_{m_{1}}H_{1}(\xi,m_{t},Du_{1}(\xi,t))(x)m_{1}(\xi,t)d\xi \, ,
	\notag
\end{align}
with terminal condition
\begin{align*}
	u_{1}(x,T)
	&=
	h_{1}(x,m_{T})+\int_{\RR^{n}}\partial_{m_{1}}h_{1}(\xi,m_{T})(x)m_{1}(\xi,T)d\xi  \, .
\end{align*}
A similar equation holds for $u_2$ and the FP equations remain~\eqref{eq:NMFC-sys-FP1}--\eqref{eq:NMFC-sys-FP2}.
Notice the difference between~\eqref{eq:NMFC-simplecase2-sys-HJB11} and the corresponding equation for CMFC, namely~\eqref{eq:CMFC-HJBFPsys-HJB}.

% subsec-subsec-subsec-subsec-subsec-subsec-subsec-subsec-subsec %
% subsec-subsec-subsec-subsec-subsec-subsec-subsec-subsec-subsec %
\subsection{Aversion in crowd motion}

We briefly revisit examples of aversion in crowd motion, in the framework of this paper. For the FP equations, we let
$$
	 A_i^* m^i = -\frac{\sigma^2}{2} \Delta m^i, \qquad g_i(x,m,v_i) = v_i,
$$
where $\sigma>0$ is a constant. In other words, the drift is the control, and the diffusion is a standard diffusion with constant volatility $\sigma$.
	In particular, $g_1$ (resp. $g_2$) does not depend upon $m_2$ (resp. $m_1$).

\paragraph{Model of Lachapelle and Wolfram~\cite{LachapelleWolfram-2011-MFG-congestion-aversion}. }
In~\cite{LachapelleWolfram-2011-MFG-congestion-aversion}, the authors considered costs of the form
\begin{align*}
	&f_i(x,m,v_i) = \chi(v_i) + \varphi_{i,\lambda}(m(x)) \, , 
	\qquad 
	h_i(x,m) = \psi_i(x) \, ,
	\\
	&
	\chi(v) = \frac{|v|^2}{2} \, , 
	\quad \varphi_{1,\lambda}(\mu_1,\mu_2) = \mu_1 + \lambda \mu_{2}, \quad \varphi_{2,\lambda}(\mu_1,\mu_2) = \mu_2 + \lambda \mu_{1} \, , 
\end{align*}
where $\lambda > 0$ is a constant and $\psi_i: \RR^n \to \RR$. In $f_i$, the terms $\varphi_{i,\lambda}$ model aversion (of an agent towards its own population or the other population). Notice that these costs are \emph{local} in $m$, in the sense that they depend only on the value of the density at the point $x$ under consideration.
The FP equation for $m_i$ becomes:
$$
	\frac{\partial m_i}{\partial t}(x,t)
	-\frac{\sigma^2}{2}\Delta m_i(x,t)
	+\diver_x \, (v_i(x,t)m_i(x,t))
	=0 \, .
$$
The Hamiltonians defined by~\eqref{eq:defH-CMFC} are
\begin{equation*}
	H_{1}(x,m,q_{1})
	=
	- \frac{1}{2} |q_1|^2 + m_1(x) + \lambda m_2(x) \, ,
	\qquad
	H_{2}(x,m,q_{2})
	=
	- \frac{1}{2} |q_2|^2 + m_2(x) + \lambda m_1(x) \, .
%	\label{eq:defH-CMFC}
\end{equation*}

For the CMFC model, the adjoint equations are
\begin{align*}
	& -\frac{\partial u_1}{\partial t}(x,t) - \frac{\sigma^2}{2} \Delta u_1(x,t) 
	 = - \frac{|\grad u_1(x)|^2}{2} + 2 \left[ m_1(x) + \lambda m_2(x) \right]
	\\
	& -\frac{\partial u_2}{\partial t}(x,t) - \frac{\sigma^2}{2} \Delta u_2(x,t) 
	 = - \frac{|\grad u_2(x)|^2}{2} + 2 \left[ m_2(x) + \lambda m_1(x) \right] \, ,
\end{align*}
whereas for the NMFC model,  if the controls are allowed to depend only on $x$ (as in~\cite[Proposition 4.1]{LachapelleWolfram-2011-MFG-congestion-aversion}; see also~\cite[Chapter 4, Proposition 4.2.1]{Lachapelle-2010-phdthesis} and its proof), the adjoint equations are
\begin{align*}
	& -\frac{\partial u_1}{\partial t}(x,t) - \frac{\sigma^2}{2} \Delta u_1(x,t) 
	 = - \frac{|\grad u_1(x)|^2}{2} + 2 m_1(x) + \lambda m_2(x)
	\\
	& -\frac{\partial u_2}{\partial t}(x,t) - \frac{\sigma^2}{2} \Delta u_2(x,t) 
	 = - \frac{|\grad u_2(x)|^2}{2} + 2 m_2(x) + \lambda m_1(x) \, .
\end{align*}
We recover the equations found in~\cite[Proposition 4.1]{LachapelleWolfram-2011-MFG-congestion-aversion} and, as noticed by Lachapelle and Wolfram, the adjoint equations for CMFC and NMFC are equivalent up to multiplying $\lambda$ by a constant.

For both CMFG and NMFG, the system of adjoint equations is given by
\begin{align*}
	& -\frac{\partial u_1}{\partial t}(x,t) - \frac{\sigma^2}{2} \Delta u_1(x,t) 
	 = - \frac{|\grad u_1(x)|^2}{2} + m_1(x) + \lambda m_2(x)
	\\
	& -\frac{\partial u_2}{\partial t}(x,t) - \frac{\sigma^2}{2} \Delta u_2(x,t) 
	 = - \frac{|\grad u_2(x)|^2}{2} + m_2(x) + \lambda m_1(x) \, .
\end{align*}

If the controls are allowed to depend on both $x$ and $m$, the necessary optimality conditions for NMFC are expressed in terms of master equations: the system~\eqref{eq:2-26}--\eqref{eq:2-29} rewrites, in this setting:
\begin{align*}
	&-\dfrac{\partial U_{1}}{\partial t}(x,m,t) 
	 - \frac{\sigma^2}{2} \Delta_x U_{1}(x,m,t) 
%	\label{eq:2-26}
	\\
	= \, 
	& \sum_{j=1}^2\int_{\RR^{n}}  \frac{\sigma^2}{2} \Delta_\xi \partial_{m_{j}}U_{1}(\xi,m,t)(x)m_{j}(\xi)d\xi   
%	\notag
%	\\
%	& 
	- \sum_{j=1}^2 \int_{\RR^{n}}\partial_{m_{j}}D_{\xi}U_{1}(\xi,m,t)(x). D_{\xi}U_{j}(\xi,m,t) m_{j}(\xi)d\xi
	\notag
	\\
	& -\frac{1}{2} |D_{x}U_{1}(x,m,t)|^2 + 2m_1(x) + \lambda m_2(x) 
%	\notag
%	\\
%	& 
	-\int_{\RR^{n}}DU_{-1}(\xi,m,t).
		\partial_{m_{1}}D_{\xi}U_{2}(\xi,m,t)(x)
	m_{2}(\xi)d\xi \, ,
	\notag
\end{align*}
and
\begin{align*}
	&-\dfrac{\partial U_{-1}}{\partial t}(x,m,t)  - \frac{\sigma^2}{2} \Delta_x U_{-1}(x,m,t) 
%	\label{eq:2-27}
	\\
	= \, 
	& \sum_{j=1}^2\int_{\RR^{n}}  \frac{\sigma^2}{2} \Delta_\xi \partial_{m_{j}}U_{-1}(\xi,m,t)(x)m_{j}(\xi)d\xi 
%	\notag
%	\\
%	& 
	-\sum_{j=1}^2 \int_{\RR^{n}}\partial_{m_{j}}D_{\xi}U_{-1}(\xi,m,t)(x). D_{\xi}U_{j}(\xi,m,t) m_{j}(\xi)d\xi
	\notag
	\\
	& - D_{x}U_{-1}(x,m,t) . D_{x}U_{2}(x,m,t) + \lambda m_{1}(x)
%	\notag
%	\\
%	& 
	-\int_{\RR^{n}}DU_{-1}(\xi,m,t).
		\partial_{m_{2}}D_{\xi}U_{2}(\xi,m,t)(x)
	m_{2}(\xi)d\xi \, ,
	\notag
\end{align*}
complemented with terminal conditions 
and analogous equations for $U_{2},U_{-2}$.

\paragraph{Model of Aurell and Djehiche~\cite{MR3763083}. }

In~\cite{MR3763083}, the authors considered a variant of the above model with \emph{non-local} running cost of the form
\begin{align*}
	&f_i(x,m,v_i) = \chi(v_i) + \Phi_{i,\Lambda}(x,m) \, , 
	\\
	&
	\chi(v) = \frac{|v|^2}{2}, 
	\quad \Phi_{i,\Lambda}(x,m) = \sum_{k=1}^2 \Lambda_{i,k} \, \phi * \mu_k(x) \, , 
\end{align*}
where $\Lambda_{i,k} \geq 0$ are constants, $*$ denotes the convolution and $\phi: \RR^n \to \RR$ is a smooth function such as
\begin{equation}
\label{eq:example-AD-phi}
	\phi(x) = \gamma_{\delta} * \mathbb{I}_{B_r}(x) \, , \qquad \gamma_{\delta}(x) = \gamma(x/\delta)/\delta \, ,
\end{equation}
where $\delta >0$, $\gamma$ is a mollifier, and $\mathbb{I}_{B_r}$ is the indicator function of the ball with radius $r$ centered at $0$ normalized by the volume of this ball. Here, $\Phi_{i,\Lambda}$  models aversion.
For the sake of comparison, let us consider the case where the controls are functions of $x$ and do not depend upon $m$ (see~\cite{MR3763083}, Assumption (C4) p. 444). One can check that this example also falls in the second special case of section~\ref{sec:NMFC:specialcases}.

For the CMFC model, the adjoint equations are
\begin{align*}
	& -\frac{\partial u_1}{\partial t}(x,t) - \frac{\sigma^2}{2} \Delta u_1(x,t) 
	 = - \frac{|\grad u_1(x)|^2}{2} 
	 + \sum_{k=1}^2 \Lambda_{1,k} \, \phi * m_{k,t}(x) + \sum_{j=1}^2 \Lambda_{j,1} \, \overline \phi * m_{j,t}(x)
	\\
	& -\frac{\partial u_2}{\partial t}(x,t) - \frac{\sigma^2}{2} \Delta u_2(x,t) 
	 = - \frac{|\grad u_2(x)|^2}{2} 
	 + \sum_{k=1}^2 \Lambda_{2,k} \, \phi * m_{k,t}(x) + \sum_{j=1}^2 \Lambda_{j,2} \, \overline \phi * m_{j,t}(x) \, ,
\end{align*}
whereas for the NMFC model, the adjoint equations are
\begin{align*}
	& -\frac{\partial u_1}{\partial t}(x,t) - \frac{\sigma^2}{2} \Delta u_1(x,t) 
	 = - \frac{|\grad u_1(x)|^2}{2} 
	 + \sum_{k=1}^2 \Lambda_{1,k} \, \phi * m_{k,t}(x) + \Lambda_{1,1} \, \overline \phi * m_{1,t}(x)
	\\
	& -\frac{\partial u_2}{\partial t}(x,t) - \frac{\sigma^2}{2} \Delta u_2(x,t) 
	 = - \frac{|\grad u_2(x)|^2}{2} 
	 + \sum_{k=1}^2 \Lambda_{2,k} \, \phi * m_{k,t}(x) + \Lambda_{2,2} \, \overline \phi * m_{2,t}(x) \, ,
\end{align*}
where $\overline \phi (x) = \phi(-x)$. Similarly to the example of~\cite{LachapelleWolfram-2011-MFG-congestion-aversion} and as noticed in~\cite{MR3763083}, the NMFC system can be seen as a CMFC system (up to changing the coefficients $\Lambda_{i,k}$ by a multiplicative constant) if one assumes that $\Lambda_{1,2}=\Lambda_{2,1}$ and $\phi$ is even (which is the case with~\eqref{eq:example-AD-phi} for example).

For both CMFG and NMFG, the system of adjoint equations is given by
\begin{align*}
	& -\frac{\partial u_1}{\partial t}(x,t) - \frac{\sigma^2}{2} \Delta u_1(x,t) 
	 = - \frac{|\grad u_1(x)|^2}{2} 
	 + \sum_{k=1}^2 \Lambda_{1,k} \, \phi * m_{k,t}(x)
	\\
	& -\frac{\partial u_2}{\partial t}(x,t) - \frac{\sigma^2}{2} \Delta u_2(x,t) 
	 = - \frac{|\grad u_2(x)|^2}{2} 
	 + \sum_{k=1}^2 \Lambda_{2,k} \, \phi * m_{k,t}(x) \, .
\end{align*}

\begin{remark}
	Here again, if the controls were allowed to depend on both $x$ and $m$, the necessary optimality conditions for NMFC would be expressed in terms of the master equations~\eqref{eq:2-26}--\eqref{eq:2-29}.
\end{remark}

%}

% subsec-subsec-subsec-subsec-subsec-subsec-subsec-subsec-subsec %
% subsec-subsec-subsec-subsec-subsec-subsec-subsec-subsec-subsec %
\subsection{Linear-quadratic models}

% secsecsecsecsecsecsecsecsecsecsecsecsecsecsecsecsecsecsecsecsecsecsecsec %
% secsecsecsecsecsecsecsecsecsecsecsecsecsecsecsecsecsecsecsecsecsecsecsec %
% secsecsecsecsecsecsecsecsecsecsecsecsecsecsecsecsecsecsecsecsecsecsecsec %
%\section{Linear-quadratic models}
%\label{secLQM}
%\setcounter{equation}{0}
%\setcounter{theorem}{0}

In this section, we consider the linear-quadratic setting. In the case of a single population, we refer the reader to, e.g., the papers~\cite{MR2352434,MR2945936,MR3489817}, the monographs~\cite{MR3134900} (Chapter 6) and~\cite{MR3752669} (Section 3.5), as well as the references therein. In the case of several populations, a model with a finite number of agents has been studied in~\cite{7403050}. Here we focus on the limit mean field models, in the different cases of interactions introduced above.

In the sequel, for $m \in L^2(\RR^n)^2$, we denote by $\overline m_k = \int x m_k(x) dx$ the first moment of $m_k$. We will denote by  $\dot{\varphi}$ the time derivative of a function $\varphi$. 

% subsec-subsec-subsec-subsec-subsec-subsec-subsec-subsec-subsec %
% subsec-subsec-subsec-subsec-subsec-subsec-subsec-subsec-subsec %
%\subsubsection{Mean field type control problems} 

We consider, for $x \in \RR^{n}, m=(m_1,m_2) \in L^2(\RR^n)^2, v_i \in \RR^d$, drift and cost functions of the form:
\begin{align}
	&g_i(x,m,v_i) = \bfA^i x + \sum_{j=1}^2 \obfA^i_j \overline m_j + \bfB^i v_i
	\label{eq:def-g-LQMFC}
	\\
	&f_i(x,m,v_i) = \frac{1}{2}\left[ x^* \bfQ^i x + (v_i)^* \bfR^i v_i + \sum_{j=1}^2 \left(x - \bfS^i_j \overline m_j\right)^* \obfQ^i_j \left(x - \bfS^i_j \overline m_j\right) \right] 
	\label{eq:def-r-LQMFC}
	\\
	&h_i(x,m) =  \frac{1}{2}\left[ x^* \bfQ^i_T x + \sum_{j=1}^2 \left(x - \bfS^i_{T,j} \overline m_j\right)^* \obfQ^i_{T,j} \left(x - \bfS^i_{T,j} \overline m_j\right) \right],
	\label{eq:def-h-LQMFC}
\end{align}
where $\bfA^i, \obfA^i_j,$ and $\bfB^i$ are bounded deterministic matrix-valued functions in time of suitable sizes, $\bfS^i_j$ (respectively $\bfQ^i, \obfQ^i_j$ and $\bfR^i$) are bounded deterministic (respectively, non-negative and positive definite) matrix-valued functions in time of suitable sizes. Except for $t=T$, we omit the dependence on time to alleviate the notations. Moreover $M^*$ denotes the transpose of a matrix $M$.

In this setting the Fokker-Planck equation~\eqref{FPCeqGen} rewrites
\beq
\label{FPCeqLQ}
\begin{split}
	\frac{\partial m_i}{\partial t}(x,t)
	+A_i^*m_i(x,t)
	+\diver_x \, \left(\left[\bfA^i x + \sum_{j=1}^2 \obfA^i_j \overline m_{j,t} + \bfB^i v_i(x,t)\right]m_i(x,t)\right)
	&=0 \, .
\end{split}
\eeq

\paragraph{Common mean field control problem. } We first investigate the CMFC problem. 
We look for adjoint states of the form 
\begin{equation}
\label{eq:LQCMFC-u}
	u_i(x,t) = \frac{1}{2} x^* P^i_t x + x^* \nu^i_t + \tau^i_t.
\end{equation}
We have $D u_i(x,t) = P^i_t x + \nu^i_t$.

The Hamiltonian~\eqref{eq:defH-CMFC} rewrites, for $x \in \RR^{n}, m=(m_1,m_2)^* \in L^2(\RR^n)^2, q_i \in \RR^n$, as
\begin{align}
	H_i(x,m,q_i)
	&=  
	q_i^* \left[\bfA^i x 
	+ \sum_{j=1}^2 \obfA^i_j \overline m_j \right]
	- \frac{1}{2} q_i^* \bfB^i (\bfR^i)^{-1} (\bfB^i)^* q_i
	\notag
	\\
	& \qquad 
	+ \frac{1}{2} \left[ x^* \bfQ^i x + \sum_{j=1}^2 \left(x - \bfS^i_j \overline m_j\right)^* \obfQ^i_j \left(x - \bfS^i_j \overline m_j\right) \right],
	\label{eq:Hk-CMFC-LQ}
\end{align}
since $\hat{\bfv}_i(x,m,q_i) = - (\bfR^i)^{-1} (\bfB^i)^* q_i$.

Using~\eqref{FPCeqLQ} and
taking into account the expression of $\hat{\bfv}$, the first moments $\overline m_i$ of $m_i$ solve the system of ODEs
\begin{equation}
\label{eq:LQCMFC-eqm1}
	\dot{\overline m}_{i,t}
	= \bfA^i {\overline m}_{i,t} + \sum_{j=1}^2 \obfA^i_j \overline m_{j,t} - \bfB^i (\bfR^i)^{-1} (\bfB^i)^* P^i_t \overline m_{i,t} - \bfB^i (\bfR^i)^{-1} (\bfB^i)^* \nu^i_t ,
\end{equation}
and the initial condition ${\overline m}_i(0) = \int x \rho_{i,0}(x) dx$.

We see that the integral term in~\eqref{eq:CMFC-HJBFPsys-HJB} rewrites
\begin{align}
\label{eq:LQCMFC-integral}
	\int \sum_{j=1}^2  \partial_{m_i} H_j(\xi,m, D u_j(\xi,t))(x) m_j(\xi,t) d\xi
	&=  
	\sum_{j=1}^2 \left[\left(P^j_{t} \overline m_{j,t} + \nu^j_{t} \right)^* \obfA^j_i - \left(\overline m_{j,t} - \bfS^j_i \overline m_{i,t} \right)^* \obfQ^j_i \bfS^j_i\right] x.
\end{align}

Replacing $u_i$ by its expression~\eqref{eq:LQCMFC-u} in the adjoint equation~\eqref{eq:CMFC-HJBFPsys-HJB}, we obtain that $(P^1,P^2),$ $ (\nu^1,\nu^2)$ and $(\tau^1,\tau^2)$ solve the following system of ODEs, which is coupled with the equations on the first moments~\eqref{eq:LQCMFC-eqm1},
\begin{align*}
	&\dot P^i_t + (P^i_t)^* \bfA^i + (\bfA^i)^* P^i_t - (P^i_t)^* \bfB^i (\bfR^i)^{-1} (\bfB^i)^* P^i_t + \bfQ^i + \sum_{j=1}^2 \obfQ^i_j = 0
	\\
	&-\dot \nu^i_t = \left[ (\bfA^i)^* - (P^i_t)^* \bfB^i (\bfR^i)^{-1} (\bfB^i)^* \right] \nu^i_t + \frac{1}{2} \sum_{j=1}^2 \left[ (P^i_t)^* \obfA^i_j + (\obfA^i_j)^* P^i_t - \obfQ^i_j \bfS^i_j - (\bfS^i_j)^*\obfQ^i_j  \right] \overline m_{j,t}
	\\
	& \qquad\qquad
	+ \sum_{k=1}^2 \left[ (\obfA^k_i)^* \left(P^k_t \overline m_{k,t} + \nu^k_t \right) - (\bfS^k_i)^* \obfQ^k_i  \left(\overline m_{k,t} - \bfS^k_i \overline m_{i,t}\right) \right]
	\\
	&-\dot \tau^i_t = \Tr\; a^i P^i_t 
	+ \sum_{j=1}^2 (\nu^i_t)^* \obfA^i_j \overline m_{j,t} 
	-  \frac{1}{2} (\nu^i_t)^* \bfB^i (\bfR^i)^{-1} (\bfB^i)^* \nu^i_t 
	+ \frac{1}{2} \sum_{j=1}^2 \left(\bfS^i_j\overline m_{j,t}\right)^* \overline Q^i_j \left(\bfS^i_j\overline m_{j,t}\right) \, ,
\end{align*}
with the terminal conditions
\begin{align*}
	&P^i_T = \bfQ^i_T + \sum_{j=1}^2 \obfQ^i_{T,j}
	\\
	&\nu^i_T = - \sum_{j=1}^2 \obfQ^i_{T,j} \bfS^i_{T,j} \overline m_{j,T}
	- \sum_{k=1}^2 ( \bfS^k_{T,i})^* \obfQ^k_{T,i} \left(\overline m_{k,T} - \bfS^k_{T,i} \overline m_{i,T} \right)
	\\
	&\tau^i_T = \frac{1}{2} \sum_{j=1}^2 \left(\bfS^i_{T,j} \overline m_{j,T}\right)^* \obfQ^i_{T,j} \bfS^i_{T,j} \overline m_{j,T} \, .
\end{align*}

We can relate this system of ODEs to a Riccati equation as follows. Let us look for $K^i$ such that
\begin{equation}
\label{eq:LQCMFC-Riccati}
	K^i_t \overline m_{i,t} = \int_{\RR^n} D u_i(x,t) m_i(x,t) dx . %\EE \left[\grad u^i(X^i_t,t)\right].
\end{equation}
Taking the derivative (with respect to time) on both sides of the above equality, using that $\nu^i = (K^i-P^i) \overline m_i$, and using the equations for $\dot P^i_t$ and $\dot \nu^i_t$ yields
\begin{align*}
		\dot K^i_t \overline m_{i,t}
	%\\
	=
	\, &\left[ K^i_t \bfB^i (\bfR^i)^{-1} (\bfB^i)^* K^i_t - K^i_t \left(\bfA^i + \obfA^i_i \right) - \left(\bfA^i + \obfA^i_i \right)^*K^i_t  - \bfQ^i - \sum_{j=1}^2 \obfQ^i_j \right.
	\\
	&\qquad \left. +\,  \obfQ^i_i \bfS^i_i + (\bfS^i_i)^*\obfQ^i_i - \sum_{j=1}^2 (\bfS^j_i)^* \obfQ^j_i \bfS^j_i \right] \overline m_{i,t}
		\\
		& \quad +\left[ - K^i_t \obfA^i_{-i} - \obfA^{-i}_{i} K^{-i}_t + (\bfS^{-i}_i)^* \obfQ^{-i}_i + \obfQ^i_{-i} \bfS^i_{-i} \right] \overline m_{-i}(t).
\end{align*}
This system of equations can be synthetically written under the following form, which turns out to be a \emph{symmetric} Riccati equation
$$
	\dot K_t = K_t \bfB K_t - (\bfA  + \obfA)^* K_t - K_t (\bfA  + \obfA) + \bfG, \qquad K_T = \bfG_T,
$$
where
\begin{equation}
\label{eq:LQCMFC-defBAAbar}
	K = \begin{pmatrix}
		K^1 & 0
		\\
		0 & K^2
		\end{pmatrix} \,,
%	\quad
	\bfB = \begin{pmatrix}
		\bfB^1 (\bfR^1)^{-1} (\bfB^1)^* & 0
		\\
		0 & \bfB^2 (\bfR^2)^{-1} (\bfB^2)^*
		\end{pmatrix} \,,
%	\quad
	\bfA = \begin{pmatrix}
		\bfA^1 & 0
		\\
		0 & \bfA^2
		\end{pmatrix} \,,
%	\quad
	\obfA = \begin{pmatrix}
		\obfA^1_1 & \obfA^1_2
		\\
		\obfA^2_1 & \obfA^2_2
		\end{pmatrix} \,,
\end{equation}
and
\begin{equation*}
	\bfG = \begin{pmatrix}
		\bfG^1_1 & \bfG^2_1
		\\
		\bfG^1_2 & \bfG^2_2
		\end{pmatrix}
\end{equation*}
with
\begin{align*}
	\bfG^i_i &= - \bfQ^i - \sum_{j=1}^2 \obfQ^i_j + \obfQ^i_i \bfS^i_i + ( \bfS^i_i )^* \obfQ^i_i - \sum_{j=1}^2 (\bfS^j_i)^* \obfQ^j_i \bfS^j_i 
	\\
	\bfG^{-i}_i &= (\bfS^{-i}_i)^* \obfQ^{-i}_i + \obfQ^i_{-i} \bfS^i_{-i} \, .
\end{align*}
and, for the final condition, $\bfG_T$ is defined similarly.

\paragraph{Mean field game. }
We focus on the system~\eqref{eq:CMFG-HJBFPsys-HJB}--\eqref{eq:CMFG-HJBFPsys-FP}. We look for $u_i$ under the form
\begin{equation}
\label{eq:LQCMFG-u}
	u_i(x,t) = \frac{1}{2} x^* P^i_t x + x^* \nu^i_t + \tau^i_t.
\end{equation}
Following the same approach as above for the CMFC problem and noting that the integral term~\eqref{eq:LQCMFC-integral} does not appear in the adjoint equation~\eqref{eq:CMFG-HJBFPsys-HJB} for CMFG, we obtain that $\overline m$ satisfies~\eqref{eq:LQCMFC-eqm1} and $P^i, \nu^i, \tau^i$ satisfy the system of ODEs
\begin{align*}
	&\dot P^i_t + (P^i_t)^* \bfA^i + (\bfA^i)^* P^i_t - (P^i_t)^* \bfB^i (\bfR^i)^{-1} (\bfB^i)^* P^i_t + \bfQ^i + \sum_{j=1}^2 \obfQ^i_j = 0
	\\
	&-\dot \nu^i_t = \left[ (\bfA^i)^* - (P^i_t)^* \bfB^i (\bfR^i)^{-1} (\bfB^i)^* \right] \nu^i_t + \frac{1}{2} \sum_{j=1}^2 \left[ (P^i_t)^* \obfA^i_j + (\obfA^i_j)^* P^i_t - \obfQ^i_j \bfS^i_j - (\bfS^i_j)^*\obfQ^i_j  \right] \overline m_{j,t}
	\\
	&-\dot \tau^i_t = \Tr\; a^i P^i_t 
	+ \sum_{j=1}^2 (\nu^i_t)^* \obfA^i_j \overline m_{j,t} 
	- \frac{1}{2} (\nu^i_t)^* \bfB^i (\bfR^i)^{-1} (\bfB^i)^* \nu^i_t 
	+ \frac{1}{2} \sum_{j=1}^2 \left(\bfS^i_j\overline m_{j,t}\right)^* \overline Q^i_j \left(\bfS^i_j\overline m_{j,t}\right),
\end{align*}
with the terminal conditions
\begin{align*}
	&P^i_T = \bfQ^i_T + \sum_{j=1}^2 \obfQ^i_{T,j}
	\\
	&\nu^i_T = - \sum_{j=1}^2 \obfQ^i_{T,j} \bfS^i_{T,j} \overline m^j_T
	\\
	&\tau^i_T = \frac{1}{2} \sum_{j=1}^2 \left(\bfS^i_{T,j} \overline m^j_T\right)^* \obfQ^i_{T,j} \bfS^i_{T,j} \overline m^j_T.
\end{align*}

Here again, we can relate this system of ODEs to a Riccati equation as follows. Let us look for $K^i$ such that
$$
	K^i_t \overline m_{i,t} = \int_{\RR^n} D u_i(x,t) m_i(x,t) dx . %\EE \left[\grad u^i(X^i_t,t)\right].
$$
Taking the derivative (with respect to time) on both sides of the above equality and using the equations for $P^i_t$ and $\nu^i_t$ yields
\begin{align*}
	\dot K^i_t \overline m_{i,t} = 
	\, &\left[ K^i_t \bfB^i (\bfR^i)^{-1} (\bfB^i)^* K^i_t - K^i_t \bfA^i - (\bfA^i)^* K^i_t - K^i_t \obfA^i_i  - \bfQ^i - \sum_{j=1}^2 \obfQ^i_j +  \obfQ^i_i \bfS^i_i \right] \overline m_{i,t}
		\\
		& -\left[ K^i_t \obfA^i_{-i}  + \obfQ^i_{-i} \bfS^i_{-i} \right] \overline m^{-i}_t.
\end{align*}
This system of equations can be written under the following compact form which turns out to be a \emph{non-symmetric} Riccati equation
\begin{equation}
\label{eq:LQCMFG-Riccati}
	\dot K_t = K_t \bfB K_t - \bfA^* K_t - K_t \bfA  - K_t \obfA + \tilde\bfG, \qquad K_T = \tilde \bfG_T,
\end{equation}
where $\bfB, \bfA,$ and $\obfA$ are defined by~\eqref{eq:LQCMFC-defBAAbar}, and
\begin{align*}
	\tilde \bfG = \begin{pmatrix}
		\tilde \bfG^1_1 & \tilde \bfG^2_1
		\\
		\tilde \bfG^1_2 & \tilde \bfG^2_2
		\end{pmatrix}
\end{align*}
with
\begin{align*}
	\tilde \bfG^i_i &= - \bfQ^i - \sum_{j=1}^2 \obfQ^i_j + \frac{1}{2}\left( \obfQ^i_i \bfS^i_i  + (\bfS^i_i)^* \obfQ^i_i \right)
	\\
	\tilde \bfG^{-i}_i &=  \frac{1}{2}\left( \obfQ^i_{-i} \bfS^i_{-i}  + (\bfS^i_{-i})^* \obfQ^i_{-i}  \right)\,, 
\end{align*}
and, for the final condition, $\tilde{\bfG}_T$ is defined similarly.

\begin{remark}
	In particular, by comparing~\eqref{eq:LQCMFC-Riccati} and~\eqref{eq:LQCMFG-Riccati}, one can see that linear-quadratic CMFC and CMFG are in general different. One can also check that the NMFC problem provides yet a different system of ODEs since this LQ model does not fall in the second special case of section~\ref{sec:NMFC:specialcases}. Details are omitted here and a study of this LQ setting with a comparison of all the ODE systems will be done elsewhere.
\end{remark}

\noindent \textbf{Acknowledgements. } 
We thank Alexander Aurell, Boualem Djehiche, Marie-Therese Wolfram and an anonymous referee for their helpful comments. Most of this work was done while the second and third authors were at NYU-Shanghai, T.H. as a visiting assistant professor and M.L. as a global postdoctoral fellow. The support of the NYU-ECNU Institute of Mathematical Sciences at NYU Shanghai is gratefully acknowledged.

\bigskip
\bigskip
{\small
\bibliographystyle{plain}
	\bibliography{2pop-bib}
}

\end{document}